\documentclass{gtart_h}  

\def\ifplaintex{\expandafter\ifx\csname documentclass\endcsname\relax}

\ifplaintex 
\else
\fi

\expandafter\ifx\csname beginpicture\endcsname\relax
\expandafter\ifx\csname documentclass\endcsname\relax
\input pictex \else
\input prepictex \input pictex \input postpictex \fi\fi

\def\gt{{\mathsurround=0pt\it $\cal G\mskip-2mu$eometry \&\ 
$\cal T\!\!$opology}}        %  journal title in recommended style

\def\gtp{{\mathsurround=0pt\it $\cal G\mskip-2mu$eometry \&\ 
$\cal T\!\!$opology $\cal P\!$ublications}}  % GT publications

\def\lognumber#1{\def\thelognumber{#1}}
\def\volumenumber#1{\def\thevolumenumber{#1}}
\def\papernumber#1{\def\thepapernumber{#1}}
\def\volumeyear#1{\def\thevolumeyear{#1}}

\def\pagenumbers#1#2{\def\startpage{#1}\def\finishpage{#2}}
\def\published#1{\def\publishdate{#1}}
\def\proposed#1{\def\theproposer{#1}}
\def\seconded#1{\def\theseconders{#1}}
\def\received#1{\def\receiveddate{#1}}

\def\accepted#1{\def\accepteddate{#1}}
\def\asciititle#1{\def\theasciititle{#1}}
\def\covertitle#1{\def\thecovertitle{#1}}

\def\asciiaddress#1{\def\theasciiaddress{#1}}

\long\def\asciiabstract#1{\long\def\theasciiabstract{#1}}
\def\asciikeywords#1{\def\theasciikeywords{#1}}

\let\\\par\let\thelognumber\relax
\let\thevolumenumber\relax\let\thepapernumber\relax
\let\thevolumeyear\relax\let\thesamplenumber\relax\let\startpage\relax
\let\finishpage\relax\let\publishdate\relax\let\receiveddate\relax
\let\reviseddate\relax\let\accepteddate\relax\let\theasciititle\relax
\let\thecovertitle\relax\let\theasciiauthors\relax\let\theasciiaddress\relax
\let\theasciiabstract\relax\let\theasciikeywords\relax
\let\theasciiemail\relax\let\theshortauthors\relax\let\theshorttitle\relax

\long\def\maketitlep{   % start of definition of \maketitlep

\count0=\startpage

\gt\hfill      %   Journal title (top left) 
\beginpicture
\setcoordinatesystem units <0.33truein, 0.33truein> point at 2.2 0.9
\setplotsymbol ({$\cal G$})
\plotsymbolspacing=9truept
\circulararc 315 degrees from 0 1 center at 0 0
\setplotsymbol ({$\cal T$})
\circulararc 315 degrees from 1 -1 center at 1 0
\endpicture
\break
{\small\ifx\thesamplenumber\relax % sample?  
Volume \else Sample
\fi\thevolumenumber\ (\thevolumeyear)
\startpage--\finishpage\nl
Published: \publishdate}
\vglue 0.5truein plus 0.4fil minus 0.1truein

{\parskip=0pt\leftskip 0pt plus 1fil\def\\{\par\smallskip}{\ifplaintex\large
\else\Large\fi\bf\thetitle}\par\medskip}   

\vglue 0pt plus 0.1fil 

{\parskip=0pt\leftskip 0pt plus 1fil\def\\{\par}{\sc\theauthors}
\par\medskip}

\vglue 0pt plus 0.1fil 

{\small\parskip=0pt\let\newline\\
{\leftskip 0pt plus 1fil\def\\{\par}{\sl\theaddress}\par}
\expandafter\ifx\theemail\relax    % email address?
\relax\else\vglue 5pt plus 0.02fil minus 2pt\def\\{\stdspace{\rm 
and}\stdspace} 
\cl{Email:\stdspace\tt\theemail}\fi
\ifx\theurl\relax                  % URL given?
\relax\else\vglue 5pt plus 0.02fil minus 2pt\def\\{\stdspace{\rm 
and}\stdspace}
\cl{URL:\stdspace\tt\theurl}\fi\par}

\vglue 7pt plus 0.3fil minus 3pt

{\bf Abstract}
\vglue 5pt plus 0.1fil minus 2pt

\theabstract

\vglue 7pt plus 0.3fil minus 3pt

{\bf AMS Classification numbers}\quad Primary:\quad \theprimaryclass

Secondary:\quad \thesecondaryclass

\vglue 5pt plus 0.3fil minus 2pt

{\bf Keywords:}\quad \thekeywords

\vglue 10pt plus 0.5fil minus 5pt

{\small  Proposed: \theproposer\hfill Received: \receiveddate\nl
Seconded: \theseconders\hfill 
\ifx\reviseddate\relax                         % paper revised?
Accepted: \accepteddate                        % no
\else
Revised: \reviseddate                          % yes
\fi}
\eject
}       %  end of definition of \maketitlep

\let\maketitlepage\maketitlep
\let\maketitle\maketitlepage

\font\phead=cmsl9 scaled 950
\font=cmsl9 scaled 1050
\font\pnum=cmbx10 scaled 913
\font\lnum=cmbx10 
\font\pfoot=cmsl9 scaled 950
\font=cmsl9 scaled 1050
\ifplaintex
\headline{\vbox to 0pt{\vskip -4.5mm\line{\small\phead\ifnum
\count0=\startpage ISSN 1364-0380 (on line)
1465-3060 (printed) \hfill {\pnum\folio}\else\ifodd\count0\def\\{ }% 
\ifx\theshorttitle\relax\thetitle\else\theshorttitle\fi\hfill{\pnum\folio}
\else\def\\{ and }{\pnum\folio}\hfill\ifx\theshortauthors\relax\theauthors
\else\theshortauthors\fi\fi\fi}\vss}}
\footline{\vbox to 0pt{\vglue 0mm\line{\small\pfoot\ifnum\count0=\startpage
\copyright\ \gtp\hfill\else
\gt, Volume \thevolumenumber\ (\thevolumeyear)\hfill\fi}\vss
}}
\else
\makeatletter
\def\@oddhead{{\small\lhead\ifnum\count0=\startpage ISSN 1364-0380 (on line)
1465-3060 (printed) \hfill {\lnum\number\count0}\else\ifodd\count0
\def\\{ }\ifx\theshorttitle\relax \thetitle \else\theshorttitle\fi\hfill
{\lnum\number\count0}\else\def\\{ and }{\lnum\number\count0}
\hfill\ifx\theshortauthors\relax 
\theauthors\else\theshortauthors\fi\fi\fi}}\def\@evenhead{\@oddhead}
\def\@oddfoot{\small\lfoot\ifnum\count0=\startpage\copyright\ \gtp\hfill\else
\gt, Volume \thevolumenumber\ (\thevolumeyear)\hfill\fi}
\def\@evenfoot{\@oddfoot}
\makeatother
\fi

\newwrite\gtoutfile
\long\gdef\makeheadfile{  %%% start of definition of \makeheadfile
{\def\\{, }\def\s{ }
\immediate\openout\gtoutfile head.xxx
\immediate\write\gtoutfile{Proxy-for: \ifx\theasciiauthors\relax
\theauthors\else\theasciiauthors\fi\s<\ifx\theasciiemail\relax\theemail\else\theasciiemail\fi>}
\immediate\write\gtoutfile{\noexpand\\}
\immediate\write\gtoutfile{Authors: \ifx\theasciiauthors\relax
\theauthors\else\theasciiauthors\fi}
{\def\\{ }\immediate\write\gtoutfile{Title: \ifx\theasciititle\relax
\thetitle\else\theasciititle\fi}}
\immediate\write\gtoutfile{Subj-class: GT or SG or MG etc}
\immediate\write\gtoutfile{MSC-class: \theprimaryclass\ifx\thesecondaryclass\relax\else, \thesecondaryclass\fi}
\immediate\write\gtoutfile{Journal-ref: Geom. Topol. \thevolumenumber
(\thevolumeyear) \startpage-\finishpage}
\immediate\write\gtoutfile{Comments: Published by Geometry and Topology at}
\immediate\write\gtoutfile{\s\s http://www.maths.warwick.ac.uk/gt/GTVol\thevolumenumber/paper\thepapernumber.abs.html}
\immediate\write\gtoutfile{\noexpand\\}
\immediate\write\gtoutfile{}
\ifx\theasciiabstract\relax
\immediate\write\gtoutfile{\theabstract}\else
\immediate\write\gtoutfile{\theasciiabstract}\fi
\immediate\write\gtoutfile{}
\immediate\write\gtoutfile{\noexpand\\}
\immediate\write\gtoutfile{}
\immediate\closeout\gtoutfile}}  %%% end of definition of \makeheadfile

\def\maketitlepage{\maketitlep\makeheadfile}
\let\maketitle\maketitlepage

\lognumber{535}
\received{4 December 2004}
\volumenumber{9}\papernumber{52}\volumeyear{2005}
\pagenumbers{2261}{2302}   
\accepted{1 November 2005}
\proposed{Ronald Stern}
\seconded{Peter Teichner, Ronald Fintushel}
\published{6 December 2005}

\usepackage{amsmath,amssymb,amsthm}
\usepackage{mathrsfs}
\usepackage{graphicx}
\usepackage{psfrag}
\usepackage{eufrak}

\def\psfraga <#1,#2> #3#4{%
\psfrag {#3}{\smash{\rlap{\kern #1 \raise #2\hbox{#4}}}}}

\def\figref#1{\hyperlink{#1anchor}{Figure~\ref*{#1}}}
\def\tabref#1{\hyperlink{#1anchor}{Table~\ref*{#1}}}
\def\anchor#1{\noindent\hypertarget{#1anchor}{\smash{$\phantom{99}$}}}

 \newtheorem{thm}{Theorem}[subsection]
 \newtheorem{cor}[thm]{Corollary}
 \newtheorem{lem}[thm]{Lemma}
 \newtheorem{prop}[thm]{Proposition}
 \theoremstyle{definition}
 \newtheorem{defn}[thm]{Definition}
 \theoremstyle{remark}

 \numberwithin{equation}{subsection}
\newcommand{\kommentar}[1]{}

\let \nc \newcommand
\let \rnc \renewcommand

\rnc \phi \varphi
\rnc \epsilon \varepsilon

\nc {\R}{{\bf R}} \nc {\C}{{\bf C}} \nc {\Z}{{\bf Z}}

\nc {\bd}{\begin{description}} \nc {\ed}{\end{description}} \nc
{\bi}{\begin{itemize}} \nc {\ei}{\end{itemize}} \nc
{\be}{\begin{enumerate}} \nc {\ee}{\end{enumerate}} \nc
{\bdm}{\begin{displaymath}} \nc {\edm}{\end{displaymath}} \nc
{\bea}{\begin{eqnarray*}} \nc {\eea}{\end{eqnarray*}} \nc
{\baa}{\begin{alignat*}} \nc {\eaa}{\end{alignat*}} \nc
{\bsp}{\begin{split}} \nc {\esp}{\end{split}} \nc
{\beq}{\begin{equation}} \nc {\eeq}{\end{equation}} \nc
{\btab}{\begin{tabular}} \nc {\etab}{\end{tabular}} \nc
{\ba}{\begin{array}} \nc {\ea}{\end{array}}

\newcommand{\cA}{{\mathcal A}}
\newcommand{\cF}{{\mathcal F}}
\newcommand{\cN}{{\mathscr N}}
\newcommand{\cM}{{\mathscr M}}
\newcommand{\cG}{{\mathcal G}}

\newcommand{\cL}{{\mathscr L}}

\newcommand{\cH}{{\mathcal H}}
\newcommand{\cP}{{\mathcal P}}

\newcommand{\fa}{\EuFrak{a}}
\newcommand{\fb}{\EuFrak{b}}
\newcommand{\fc}{\EuFrak{c}}
\newcommand{\fd}{\EuFrak{d}}
\newcommand{\APS}{{\mathcal P}}

\newcommand{\dover}[1]{\dot{#1}}

\nc{\LOT}[1]{L^2(\Omega^{0+1+2}(T;#1))}
\nc{\LO}[2]{L^2(\Omega^{0+1+2}(#1;#2))}
\nc{\HT}[1]{{H^{0+1+2}(T;#1)}} \nc{\cHT}[1]{{\cH^{0+1+2}(T;#1)}}
\nc{\stor}{D^2\times S^1}

\nc{\Om}{\Omega} \nc{\dist}{{\rm dist}}
\renewcommand{\deg}{{\rm deg}}
\nc{\homeo}{\approx} \nc{\im}{\text{\rm Im}} \nc{\ad}{\text{\rm
ad}} \nc{\vol}{\text{\rm vol}} \nc{\hol}{\text{\rm hol}}
\nc{\re}{\text{\rm Re}} \nc{\Id}{\text{\rm Id}}
\nc{\Mas}{\text{\rm Mas}} \nc{\SF}{\text{\rm SF}}
\rnc{\ker}{\text{\rm Ker}} \nc{\tr}{\text{\rm tr}}
\nc{\coker}{\text{\rm Coker}} \rnc{\hom}{\text{\rm Hom}}
\nc{\End}{\text{\rm End}} \nc{\Aut}{\text{\rm Aut}}
\nc{\lat}{(\frac{1}{2}\Z)^2} \nc{\hz}{\tfrac{1}{2}\Z}
\nc{\la}{\langle} \nc{\ra}{\rangle} \nc {\Ra}{\Rightarrow} \nc
{\La}{\Leftarrow} \nc {\lla}{\longleftarrow} \nc
{\nach}{\rightarrow} \nc {\equ}{\Leftrightarrow} \nc
{\gdw}{\Longleftrightarrow} \nc {\lra}{\longrightarrow} \nc
{\Lra}{\Longrightarrow} \nc {\lmt}{\longmapsto} \nc
{\quer}{\overline} \nc {\tensor}{\otimes} \nc {\Rt}{\widetilde\R}
\nc {\contract}{\lrcorner} \nc{\sgn}{\text{\rm sgn}}

\begin{document}

\title[A splitting formula for spectral flow]
 {A splitting formula for the spectral flow of\\\vglue -6pt\\the odd signature operator on $3$--manifolds\\coupled to a path of
$SU(2)$ connections}
\covertitle{A splitting formula for the spectral flow of\\the odd signature operator on $3$--manifolds\\coupled to a path of
$SU(2)$ connections}

\asciititle{A splitting formula for the spectral flow of the odd
signature operator on 3-manifolds coupled to a path of SU(2)
connections}

\author{Benjamin Himpel}

\address{
Mathematisches Institut, Universit\"at Bonn\\
Beringstr. 6, D--53115 Bonn, Germany}

\asciiaddress{
Mathematisches Institut, Universitaet Bonn\\
Beringstr. 6, D-53115 Bonn, Germany}

\email{himpel@math.uni-bonn.de}

\urladdr{http://www.math.uni-bonn.de/people/himpel/}

\begin{abstract}
We establish a splitting formula for the spectral flow of the odd
signature operator on a closed $3$--manifold $M$ coupled to a path of
$SU(2)$ connections, provided $M = S \cup X$, where $S$ is the solid
torus. It describes the spectral flow on $M$ in terms of the spectral
flow on $S$, the spectral flow on $X$ (with certain
Atiyah--Patodi--Singer boundary conditions), and two correction terms
which depend only on the endpoints.

Our result improves on other splitting theorems by removing assumptions on the
non-resonance level of the odd signature operator or the
dimension of the kernel of the tangential operator, and allows
progress towards a conjecture by Lisa Jeffrey in her work on Witten's
$3$--manifold invariants in the context of the asymptotic expansion
conjecture \cite{jeffrey}.
\end{abstract}

\asciiabstract{We establish a splitting formula for the spectral flow
of the odd signature operator on a closed 3-manifold M coupled to a
path of SU(2) connections, provided M = S cup X, where S is the solid
torus.  It describes the spectral flow on M in terms of the spectral
flow on S, the spectral flow on X (with certain Atiyah-Patodi-Singer
boundary conditions), and two correction terms which depend only on
the endpoints.  Our result improves on other splitting theorems by
removing assumptions on the non-resonance level of the odd signature
operator or the dimension of the kernel of the tangential operator,
and allows progress towards a conjecture by Lisa Jeffrey in her work
on Witten's 3-manifold invariants in the context of the asymptotic
expansion conjecture.}



\asciikeywords{Spectral flow, odd signature operator, gauge theory, 
Chern-Simons theory, Atiyah-Patodi-Singer boundary conditions, Maslov index}

\maketitle

\section{Introduction}

This article analyzes the spectral flow of the odd signature operator coupled to a path of
$SU(2)$ connections on a 3 manifold $M$ given a decomposition $M = S
\cup X$ with $S$ a solid torus.

The question is motivated by Edward Witten's description of certain
$3$--manifold invariants \cite{witten2} and Lisa Jeffrey's work
\cite{jeffrey} on the asymptotic expansion conjecture
\cite[Conjectures 7.6 and 7.7]{ohtsuki2}, particularly
\cite[Conjecture 5.8]{jeffrey} in the case of torus bundles over
$S^1$.

This work is preceded by Paul Kirk and Erik Klassen
\cite{kirk-klassen} where they treat the problem of computing the
spectral flow between irreducible flat connections under some
restrictions. In this paper these restrictions are removed.

In particular Kirk and Klassen analyzed the
spectral flow on a torus bundle
over the circle by way of the splitting theorem in \cite{cappell-lee-miller2} applied to a
decomposition of the manifold into a solid torus and its complement.
They showed that the spectral
flow between irreducible, flat connections is $0 \bmod 4$ \cite[Theorem 7.5]{kirk-klassen},
provided that the kernel of the tangential operator for
the odd signature operator on the splitting torus along the path
of connections has constant
dimension.

Kirk and Klassen proposed \cite[Appendix]{kirk-klassen} that it might
be possible to always to find a solid torus so that the kernel of
the tangential operator has constant dimension. It is shown in Proposition \ref{example}, that this is not
always the case. Thus we are forced to deal
with the case that the dimension changes.

There have been a lot of general splitting formulas
\cite{cappell-lee-miller2, bunke, daniel-kirk, kirk-lesch, nico}.
Unfortunately, for all practical purposes one has to make a couple of technical assumptions in order to
apply any of these splitting formulas, namely
\be
\item that the non-resonance level (see \cite{nico}) of the operator in question is
 zero, and
\item that the kernel of the tangential operator has constant
dimension.
\ee
The purpose of this article is to establish a splitting formula in the
case where one side is a solid torus
without making the above assumptions. We will also show how to apply
this tool to explicitly to compute
spectral flow.

The main results of this article are the following.
\bi
\item Section \ref{torusboundaryconditions} introduces a continuous family of Atiyah--Patodi--Singer boundary
  conditions $\cP^\pm(L)$ for
manifolds with torus boundary
(Definition \ref{LambdaL}, Theorem \ref{LambdaLcontinuous}).
\item In Section \ref{SFonClosedMfld} we consider a decomposition of a
  3 manifold $M = X \cup S$ with $S$ the solid torus and express the spectral flow of the twisted odd signature
operator between flat $SU(2)$ connections as the sum of the spectral
  flow on $S$, the spectral flow on $X$ (with boundary conditions
  $\cP^\pm(L)$), and two Maslov triple indices (Theorem
\ref{splittingformula}).
\item Section \ref{SFonSolidTorus} contains an explicit computation of
the spectral flow of the twisted odd signature operator between flat
connections on the solid torus with boundary conditions $\cP^\pm(L)$ (Theorem
\ref{SFSolidTorus}).
\item In Section \ref{SFonTorusBundle} we use Theorem
  \ref{splittingformula} to compute the spectral flow of the twisted odd signature operator between
irreducible, flat connections on torus bundles over $S^1$ (Theorem
\ref{SFOnTBundlesOverS1}). In particular we eliminate the technical assumption used
in \cite{kirk-klassen} on the dimension of the kernel of the twisted
de Rham operator.
\ei

Complete calculations and proofs can be found in the author's thesis \cite{himpel}.

\rk{Acknowledgments}
The results of this article were obtained during the author's Ph.D. studies at
Indiana University, Bloomington. The author is grateful for the guidance and
support of his advisor and friend Paul Kirk. The
comments of the referee were much appreciated.
The author thanks the Max--Planck--Institut in Bonn for their
hospitality and financial support while this paper was written.

\section{Preliminaries}

Familiarity with \cite{boden-herald-kirk-klassen, kirk-klassen,
  kirk-lesch, nico} is useful. However, we introduce all the necessary facts for the convenience of the reader.

\subsection{Setup}\label{conventions}

For the rest of the article let us assume the following.

\begin{figure}[ht!]\anchor{collars}\small
\begin{center}
\leavevmode
\psfrag{S}{$S$} \psfrag {X}{$X$} \psfraga <0pt,3pt> {Sigma}{$T$}
\psfrag{[-1,1]}{$[-1,1]$}
\includegraphics[scale=.4]{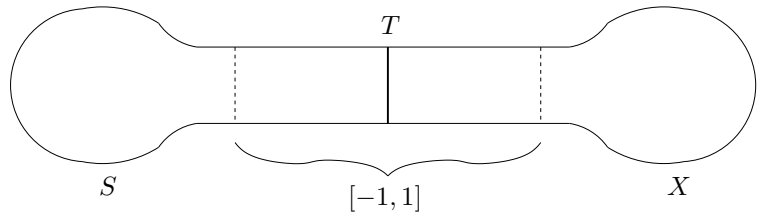}
\end{center}
\caption{\label{collars}The collar around $T$}
\end{figure}

\be
\item The orientation of the torus $T = S^1 \times S^1 = \{(e^{im},e^{il})\mid m,l \in
 [0,2\pi)\}$  is determined by $dm \wedge dl \in
\Omega^2(T)$. $T$ is given the product
metric with the standard metric on $S^1$. The fundamental group $\pi_1 T$ is the free abelian group generated by $\mu =
 \{(e^{im},1)\}$ and $\lambda = \{(1,e^{il})\}$.
\item The solid torus $S = D^2 \times S^1 = \{(n
  e^{im},e^{il}) \mid n \in [0,1], m,l \in
[0,2\pi)\}$ is
  oriented so that $dn\wedge dm \wedge dl \in
\Omega^3 (S)$ is a positive multiple of the volume form when $n >
  0$ and we have $\partial S = T$ as oriented manifolds (outward normal first convention). Consider a metric on $D^2$ such
  that a collar of $S^1 = \partial D^2$ may be
isometrically identified with $[-1,0]\times S^1$. $S$ is equipped with the
  product metric of our metrics on $D^2$ and $S^1$. Then there
  is a collar $N(T)$ of $\partial S$ which is isometric to $[-1,0]\times T$. This is consistent
  with our metric on $T$. The
  fundamental group $\pi_1 S$ is infinite cyclic generated by
  $\lambda$, and $\mu$ is trivial in $\pi_1 S$.
\item The $3$--manifold $X$ has boundary $T$ and is oriented so
  that $\partial X = - T$ as oriented manifolds. Consider a metric on $X$ such that a collar of $\partial X$ is isometric to $[0,1]\times
  T$.
\item Consider the $3$--manifold $M = S \cup_T X$ with the orientation
  and metric induced by the orientation and metric on $S$ and $X$. See
  \figref{collars}. In
  Section \ref{SFonTorusBundle}, $M$ is a torus bundle over $S^1$.
\item Let $P$ be a trivialized principal bundle with structure group
  $SU(2)$ over
  $M$ and consider its restriction to $T$, $S$ and $X$.
\ee

\subsection{Spectral flow}

Let $D_t$, $t\in [0,1]$, be a $1$--parameter family of
(possibly unbounded)
self-adjoint operators with compact resolvent, continuous in the graph
metric. Then, given $\epsilon > 0$ smaller than the modulus of the
largest negative eigenvalue of $D_0$ and $D_1$, the {\em
spectral flow} $\SF(D_t)\in \Z$ is roughly defined to be the
algebraic intersection number in $[0,1]\times \R$ of the track of
the spectrum
$$
\{(t,\lambda) \mid t\in [0,1],\lambda \in \text{Spec}(D_t)\}
$$
and the line segment from $(0,-\epsilon)$ to $(1,-\epsilon)$. The
orientations are chosen so that if $D_t$ has spectrum $\{n+t \mid
n\in Z\}$ then $\SF(D_t)=1$. For a detailed discussion on spectral flow
with rigorous definitions and proofs see
\cite{booss-lesch-phillips}.

Given a continuous family $D_t$ of self-adjoint elliptic operators on
a closed manifold $M$, one can use standard theorems
(for example in \cite{kato}) to conclude that the spectrum is discrete for each
$t$ and varies continuously in $t$.

For a compact manifold $Y$ with
boundary $\Sigma$ one needs to choose boundary conditions. In
\cite{booss-lesch-phillips} the authors treat this in more
generality, but for the results in the present article it is sufficient to confine
ourselves to Atiyah--Patodi--Singer
boundary conditions (see Section
\ref{twistedlaplacian} for a definition)
on a manifold with a collar. It is shown in
\cite{kirk-klassen2}, that by choosing a continuous family of Atiyah--Patodi--Singer
boundary conditions the spectrum of $D_t$ varies continuously in $t$
(see \cite[Theorem 4.1 and Section 5]{kirk-klassen2}).

Observe, that there are other conventions for spectral
  flow (see for example \cite{jeffrey} and
  \cite{kirk-klassen}). However, they differ only by $\dim(\ker(D_0))$ or $\dim(\ker(D_1))$.

\subsection{Connections and representations}

There are several equivalent definitions for connections, but in the
case of principal $SU(2)$ bundles over $3$--manifolds we have a
particularly nice description for connections.

We will identify the Lie group $SU(2)$ with the unit
quaternions $$\{ v \in \R \oplus \R i \oplus \R j \oplus \R k \mid |v|
= 1\}.$$ The group multiplication will be the usual multiplication in
the quaternions. Its Lie algebra $su(2)$ can be identified with the vector
space $\R i \oplus \R j \oplus
\R k = \R i \oplus \C j$ of imaginary quaternions. The adjoint action
$\ad$ of $SU(2)$ on
$su(2)$ corresponds to the conjugation of imaginary quaternions by
unit quaternions.

Consider the principal $SU(2)$ bundle $P$ over $M$. The chosen trivialization of $P$ induces an identification
\bi
\item of the group of gauge transformations with
  $\cG_M := C^\infty(M,SU(2))$, and
\item of the affine space of connections  with Lie-algebra-valued $1$--forms
$\cA_M := \Omega^1(M)\tensor su(2)$ on $M$,
where the product (or trivial) connection corresponds to
$0\in\Omega^{1}(M)\tensor su(2)$.
\ei

Then the action of a
gauge transformation $g \in C^\infty(M,SU(2))$ on an
$SU(2)$ connection $A\in\Omega^1(M)\tensor su(2)$ on $M$ is given by $
g\cdot A = g A g^{-1} + g\, dg^{-1}.$

The choice of
trivialization $P =
M \times SU(2)$ induces a trivialization of the adjoint bundle $P
\times_\ad su(2)$ associated to $P$. Thus
we can associate to a connection $A$ on $P$ a {\em covariant
derivative} on $\ad P = M \times su(2)$ by defining
\bea
d_A \co \Om^0(M) \tensor su(2) & \to & \Om^1(M) \tensor su(2)\\
a & \mapsto & d a + [A,a],
\eea
where $[\cdot,\cdot]$ means taking the Lie bracket on the coefficients
and wedging the form part. We extend the covariant derivative to
$\Om^*(M) \tensor su(2)$ by the Leibniz rule.

The {\em curvature} $F_A$ of a connection $A \in
  \Om^1(M)\tensor su(2)$ is the $2$--form in $\Om^2(M)
  \tensor su(2)$ defined by $[F_A
  ,\phi] := d^2_A \phi$, where $\phi \in \Om^0(M)\tensor su(2)$. We have
  $F_A = dA + A
\wedge A$. We call $A$ {\em flat}, if $F_A=0$. Let $\cF_M \subset \cA_M$ be
  the set of all flat connections on $M$. If $A$ is flat, then the {\em holonomy}
  $\hol(A)\co \pi_1M \to SU(2)$ is a $SU(2)$ representation of the
  fundamental group $\pi_1 M$ of $M$.

Let $\chi_M$ be the set of
conjugacy classes of $SU(2)$ representations of the fundamental group
$\pi_1(M)$. Then the holonomy map identifies $\cF_M/\cG_M$ with $\chi_M$ \cite[Proposition 2.2.3]{donaldson-kronheimer}.

\subsection{The odd signature operator and the de Rham operator}

We introduce two first order differential operators on $M$ and $T$. These depend on the orientation and
the Riemannian metric on $M$ and $T$, which we fixed in Section \ref{conventions}. We set $\Omega^{i}(M;su(2)):=
\Omega^i(M)\tensor su(2)$,
$\Omega^{0+1}(M;su(2)):=\Omega^0(M;su(2)) \oplus
\Omega^1(M;su(2))$.

The $L^2$ inner product on $\Omega^{0+1}(M;su(2))$ and
$\Omega^{0+1+2}(T;su(2))$ is given in terms of the Hodge $*$ operator by
$$\la \alpha,\beta \ra_{L^2(M)} = - \int_M \tr (\alpha \wedge *
 \beta) \quad \text{ and } \quad \la \alpha,\beta \ra_{L^2(T)} = - \int_T \tr (\alpha \wedge *
  \beta).$$
For an $SU(2)$ connection $A\in\Omega^{1}(M;su(2))$ the {\em odd signature
operator twisted by $A$} is defined to be
\begin{eqnarray*}
D_A\co \Omega^{0+1}(M;su(2)) & \to & \Omega^{0+1}(M;su(2))\\
(\alpha,\beta) & \mapsto & (d_A^*\beta,*d_A \beta + d_A \alpha).
\end{eqnarray*}
For an $SU(2)$ connection $a\in\Omega^{1}(T; su(2))$ {\em de Rham
operator twisted by $a$} is defined to be
\begin{eqnarray*}
S_a\co \Omega^{0+1+2}(T;su(2)) & \to & \Omega^{0+1+2}(T;su(2))\\
(\alpha,\beta,\gamma) & \mapsto & (*d_a\beta,-*d_a \alpha-d_a*
\gamma, d_a* \beta).
\end{eqnarray*}
In fact, $S_a$ is the {\em
tangential operator} of $D_A$ (see \cite[Lemma 2.4]{boden-herald-kirk-klassen}).
$D_A$ and $S_a$ are first-order self-adjoint elliptic differential
 operator.

The following important and well-known fact is an application of the Atiyah--Patodi--Singer index
  theorem (compare with \cite[Theorem 7.1]{kirk-klassen-ruberman}):

\begin{prop} \label{degreeshift} Suppose $M$ is a closed $3$--manifold and $g\co M \to SU(2)$
  is a gauge transformation. If $A_t$ is any path of $SU(2)$ connections on $M$ from $A_0$ to $A_1 = g
  \cdot A_0$, then
$$
\SF(D_{A_t}) = 8 \, \deg(g)
$$
\end{prop}

\subsection{Eigenspaces of the de Rham operator}\label{twistedlaplacian}

Let $a\in \cA_T$ and $\nu \ge 0$. Denote by $E_{a,\nu}$ the $\nu$--eigenspace of $S_a$ and let
\bea
P^+_{a,\nu} & := & \text{span}_{L^2} \{\psi \mid S_a \psi = \mu \psi,
\mu> \nu \} = \overline{\bigoplus_{\mu >\nu} E_{a,\mu}}^{L^2},\\
P^-_{a,\nu} & := & \text{span}_{L^2} \{\psi \mid S_a \psi = \mu \psi,
\mu< -\nu \} = \overline{\bigoplus_{\mu <-\nu} E_{a,\mu}}^{L^2},\\
E^+_{a,\nu} & := & \bigoplus_{0<\mu \le \nu} E_{a,\mu} \quad\text{and}\\
E^-_{a,\nu} & := & \bigoplus_{-\nu\le\mu <0 } E_{a,\mu}.
\eea
Abbreviate $P^\pm_a  :=  P^\pm_{a,0}$. By the spectral theorem for
self-adjoint elliptic operators we have
\begin{equation*}%\label{splittingP}
L^2(\Om^{0+1+2}(T,su(2))) = P^+_a \oplus \ker S_a \oplus P^-_a.
\end{equation*}
We define an almost complex structure on $\LO{T}{su(2)}$ by
$$J(\alpha,\beta,\gamma):=(-*\gamma,*\beta,*\alpha).$$
One can
see that $J^2=-\Id$
and that $J$ is an isometry of $\LO{T}{su(2)}$. The induced
symplectic structure $\omega(x,y):= \la x,J y \ra$ is compatible with
 $J$, which makes $\LO{T}{su(2)}$ into a {\em symplectic space with
 compatible almost complex structure}. A subspace
$\Lambda$ of $\LO{T}{su(2)}$ is {\em
  Lagrangian} if $J \Lambda = \Lambda^\perp$.

$\ker S_a $ is a finite
  dimensional symplectic subspace of $\LO{T}{su(2)}$. Furthermore, if
  $L$ is a Lagrangian subspace of $\ker S_a$, then $L\oplus  P^{\pm}_{a}$ is Lagrangian in
  $\LO{T}{su(2)}$.

A Lagrangian subspace $\cP\subset\LO{\Sigma}{su(2)}$ is called an
  {\em Atiyah--Patodi--Singer (APS) boundary condition}, if $\cP$ contains all eigenvectors of
 the tangential operator $S_a$ with sufficiently large eigenvalues.

Let $\Delta_a := d_a d_a^* + d_a^* d_a$
  be the Laplacian
  on $\Omega^{0+1+2}(T;su(2))$ twisted by $a\in
  \cA_T$. Denote the harmonic $0$, $1$ and $2$--forms
of $\Delta$ and $\Delta_a$ by $\cH^{0+1+2}(T;su(2)) := \ker \Delta$ and
$\cH_a^{0+1+2}(T; su(2)) := \ker \Delta_a$ respectively.

If $a$ is flat, we have $d_a^2 = 0$ and consequently $\Delta_a = S_a^2$. Thus a $\lambda$--eigenvector  $\phi$ of $S_a$
is a $\lambda^2$--eigenvector of $\Delta_a$. Furthermore, a direct computation shows
the following.

\begin{lem}\label{eigcorr} If $a$ is flat and $\phi$ is a $\lambda^2$--eigenvector ($\lambda > 0$)
for $\Delta_a$, then $\phi\pm\frac{1}{\lambda}S_a\phi$ is a
$\pm\lambda$--eigenvector for $S_a$. Furthermore $\ker S_{a}=\cH_a^{0+1+2}(T; su(2))$.
\end{lem}

\subsection{Cauchy data spaces and adiabatic limits}

Cauchy data spaces play an important role in this work because
of their relation to spectral flow. Liviu Nicolaescu analyzed this
relationship in \cite{nico}. His results have been extended by
\cite{daniel-kirk} and \cite{kirk-lesch}. The facts in this section
make up the main tools for this work.

We will state all results in this section in terms of the odd
signature operator on $M$,
but they apply to other self-adjoint Dirac type operators as
well.

We introduce the notation
\bea
S^R = S \cup ([0,R]\times
T) & \text{ and } & S^\infty = S \cup ([0,\infty)\times
T),\\
X^R = X \cup ([-R,0]\times
T) & \text{ and } & X^\infty = X \cup ((-\infty,0]\times
T)
\eea
for $R \ge 0$.

Let $A$ be a connection on $X$, which is in cylindrical
form in a collar of $T$, that is $A = i^*_u a$, where $i_u\co T \hookrightarrow [-1,1] \times T$ is
the inclusion at $u \in [-1,1]$ and $a\in \cA_T$.

We write the restriction of
  $\Omega^{0+1}([-1,1]\times T;su(2))$ to $T$ as
\bea r\co \Omega^{0+1}([-1,1]\times T;su(2)) & \to &
  \Omega^{0+1+2}(T;su(2))\\
(\sigma,\tau) & \mapsto & (i_0^*(\sigma),i^*_0(\tau),*i_0^*(\tau
\contract
  \frac{\partial}{\partial u})),
\eea where $\tau\contract \frac{\partial}{\partial u}$ denotes
  contraction of $\tau$ with $\frac{\partial}{\partial u}$, and $*$ is
the Hodge star on differential forms on the $2$--manifold $T$.
This also gives us a restriction map of $\Omega^{0+1}(S;su(2))$ and  $\Omega^{0+1}(X;su(2))$ to
$\Omega^{0+1+2}(T;su(2))$.
If we write $\tau = \beta + \omega\, du$, where $u$ is the
coordinate in $[0,1]$ and $\beta$ does not have a $du$ component,
then a more intuitive way to write the restriction map is
$r(\sigma, \beta + \omega du) = (\sigma |_T, \beta|_T,
*(\omega|_T))$.

The {\em Cauchy data space} of $D_A$ is $$\Lambda_{X,A}:=\Lambda_X(D_A) :=
  \overline{r(\ker D_A)}^{L^2},$$
and the {\em scattering Lagrangian} or the {\em limiting values of extended
  $L^2$ solutions} is
$$\cL_{X,A} := \text{proj}_{\ker S_a}
(  \Lambda_{X,A} \cap (P^- \cup \ker S_a)).$$
See \cite{bleecker-booss} for more
  information on Cauchy data spaces, particularly \cite[Definition 2.22]{bleecker-booss}. Note
  that we can
  extend $D_A$ to $X^R$ and that
  $\Lambda_{X,A}^R:=\Lambda_{X^R}(D_A)$ is a continuous family of
  Lagrangian subspaces by \cite[Lemma 3.2]{daniel-kirk}. Denote
  $\Lambda^\infty_{X,A}:=\lim_{R\to\infty}\Lambda^R_{X,A}$.

For $a$ flat the kernel of $S_a$ is isomorphic to the cohomology
$H^*(\partial X;su(2)_{\hol(a)})$ with values in $su(2)$ twisted by
$\hol(a)$ via the Hodge and de
Rham theorems. See for example \cite[Chapter 5]{davis-kirk} for a
  definition of the {\em cohomology of $X$ twisted by
  a representation $\rho\co \pi(X) \to \Aut(V)$}.

The following gives a way to compute the scattering
Lagrangian at a flat connection.
\begin{prop}[Corollary 8.4, \cite{kirk-lesch}]\label{scattidentif} If $A$ is flat on
  $X$, then $\cL_{X,A}$ is isomorphic to
  $\im (H^*(X;su(2)_{\hol(A)}) \to H^*(\partial X;su(2)_{\hol(A)}))$
  via the Hodge and de Rham theorems.
\end{prop}

Cauchy data spaces are complicated and we wish to relate them to
simpler Lagrangians. The situation is particularly favorable when $D_A$ has non-resonance
level $0$. Nicolaescu \cite{nico} defines $D_A$ to have {\em non-resonance level} $\nu\ge 0$, if $P^-_{a,\nu} \cap \Lambda_{X,A} =
0$.

Nicolaescu's adiabatic limit theorem describes the limit of the
Cauchy data
spaces of $D_A$ under stretching the collar.

\begin{thm}[Corollary 4.11, \cite{nico}]\label{nicoadiabaticlimit} If $D_A$ has non-resonance level $0$, then
$$
\Lambda^\infty_{X,A} = P^+_{a} \oplus \cL_{X,A}
$$
\end{thm}

We will need some facts about the $0$ non-resonance level
situation. In analogy to \cite[Proposition
  2.10]{boden-herald-kirk-klassen} we have the following from
\cite[Proposition 4.9]{atiyah-patodi-singer}.
\begin{prop}\label{sescor} Let $A$ be flat.
\be
\item We have
$$
\Lambda_{X,A} \cap P^-_a \cong \im (H^1(X,T; su(2)_{\hol(A)}) \to H^1(X;su(2)_{\hol(A)})).
$$
Thus $0$ non-resonance level is equivalent to $$\im (H^1(X,T; su(2)_{\hol(A)})
\to H^1(X;su(2)_{\hol(A)})) = 0.$$
\item Assuming $0$ non-resonance level, we get the isomorphism
$$
\Lambda_{X,A} \cap (P^-_a \oplus Q) \cong \cL_X \cap Q.
$$
\ee
\end{prop}

Kirk and Lesch \cite{kirk-lesch} give a detailed decomposition of $L^2(\Om^{0+1+2}(T;su(2)))$ in the spirit of the
Hodge decomposition as an orthogonal sum of symplectic
spaces
\begin{equation}
\label{refineddecomposition}
(P^-_{a,\nu} \oplus (P^+_{a,\nu})) \oplus
( d_a(E^+_{a,\nu})\oplus d_a^*(E^-_{a,\nu}) ) \oplus (d_a^*(E^+_{a,\nu}) \oplus
d_a(E^-_{a,\nu})) \oplus \ker S_a .
\end{equation}
When $A$ is flat, we do not only get an explicit description of the
scattering Lagrangian (Proposition \ref{scattidentif}), but the
adiabatic limit has a nice description when the the non-resonance level is not $0$.

\begin{thm}[Theorem 8.5, \cite{kirk-lesch}]\label{adiabaticlimit}
If $A$ is flat, then there exists a subspace
$$
W_a \subset d_a(E^+_{a,\nu}) \subset P^-_{a,0}
$$
isomorphic to
$$
\im\left(H^{0+1}(X,T;su(2)_{\hol(A)}) \to H^{0+1}(X;su(2)_{\hol(A)})\right)
$$
so that if $W^\perp_a$ denotes the orthogonal complement of $W_a$ in
$d_a(E_{a,\nu}^+)$, then with respect to the decomposition
(\ref{refineddecomposition}) into symplectic subspaces, the adiabatic
limit of the Cauchy data spaces decomposes as a direct sum of
Lagrangian subspaces:
$$
\Lambda^\infty_{X,A} = P^+_{a,\nu} \oplus ( W_a \oplus
J(W_a^\perp)) \oplus d_a(E^-_{a,\nu}) \oplus \cL_A
$$
where $\cL_A \subset \ker S_a \cong H^*(T;su(2)_{\hol(A)})$
denotes the scattering Lagrangian on $X$.
\end{thm}

In due course, we will need to have control over the intersection of the Cauchy data spaces
with other Lagrangian subspaces when we stretch the collar. Thus the
following results are important.

\begin{prop}[Lemma 8.10,
    \cite{kirk-lesch}]\label{intersectionstretching} Let $A$ be flat,
    and consider a Lagrangian subspace $V \subset \ker S_{a}$.
\be
\item The dimension of $\Lambda^R_{S,A} \cap
  \Lambda^R_{X,A}$ is independent of $R\in [0,\infty]$.
\item The dimension of $\Lambda^R_{S,A} \cap
  (P_a^+ \oplus V)$ is independent of $R\in [0,\infty]$.
\item The dimension of $(P_a^- \oplus V) \cap
  \Lambda^R_{X,A}$ is independent of $R\in [0,\infty]$.
\ee
\end{prop}

\subsection{Maslov index}\label{maslovindex}

Let $H$ be a symplectic Hilbert space with compatible almost complex
structure $J$.
A pair of Lagrangians $(L,M)$ in $H$ is called {\em Fredholm} if $L+M$
is closed and both $\dim(L\cap M)$ and $\text{codim}(L+M)$ are
finite. Consider a continuous path $(L_t, M_t)$ of Fredholm pairs of
Lagrangians in $H$. Continuity is measured in the gap topology. If
$L_t$ and $M_t$ are transverse at the end points, that is, intersect trivially, then the {\em Maslov
index} $\Mas(L_t,M_t)$ is roughly defined to be a count of how many times $L_t$ and $M_t$
intersect with sign and multiplicity, that is, counting the
dimension of the intersection. For a careful definition see
\cite{cappell-lee-miller, nico, daniel1}. If the
intersection of $L_t$ and $M_t$ is nontrivial at the endpoints, we
will choose a convention compatible with our spectral flow
convention (in view of Theorem \ref{nicoextension}): Given a continuous
  $1$--parameter family of Fredholm pairs of Lagrangians $(L_t,M_t)$,
  $t\in[0,1]$, choose $\epsilon > 0$ small enough so that
\be
\item $e^{s J} L_i$ is transverse to $M_i$ for $i=0,1$ and $0<s\le
  \epsilon$, and
\item $(e^{s J} L_t, M_t)$ is a Fredholm pair for all $t\in [0,1]$ and
  all $0\le s \le \epsilon$.
\ee
Then the {\em Maslov index} of the pair $(L_t,M_t)$ is the
  Maslov index of $(e^{\epsilon L}L_t,M_t)$.

A splitting theorem for spectral flow by Nicolaescu \cite{nico} in terms of a Maslov index has been
extended in \cite{daniel} to the situation when the Dirac operators are not invertible
at the endpoints. Also see \cite[Theorem 7.6]{kirk-lesch} for a
proof of the same result. The precise statement in the context of the odd signature operator and
$SU(2)$ connections on $M = S \cup_T X$ is the following.

\begin{thm}[Theorem 4.3, \cite{daniel}]\label{nicoextension} Suppose $A_t$ is a
  continuous path of $SU(2)$ connections on $M$ in cylindrical form in
  a collar of $T$. Then the path $(\Lambda_S(D_{A_t})),
  \Lambda_X(D_{A_t}))$ consists of Fredholm pairs of Lagrangians and
$$
\SF(D_{A_t}) = \Mas(\Lambda_S(D_{A_t}),
  \Lambda_X(D_{A_t})).
$$
\end{thm}

We also have a relative version of this theorem (see \cite{daniel1} and
\cite{nico}) which relates spectral flow on a manifold with boundary
with APS boundary conditions to some Maslov index. It is implied by
the results in \cite{daniel}.

\begin{thm}\label{nicorelative} If
  $A_t$ is a path of connections on $X$ in cylindrical form near
  $T$ and $\APS_t$ is a continuous family of self-adjoint APS
  boundary conditions, then the spectral flow $\SF(D_{A_t}|X;\APS_t)$
  is well defined and $\SF(D_{A_t}|X;\APS_t) = \Mas(\Lambda_X(D_{A_t}),\APS_t)$.
\end{thm}

We will also use a Maslov triple index as defined in
\cite{kirk-lesch}, which is up to normalization the same as Bunke's
Maslov triple index in \cite{bunke} and different from the Maslov
triple index $\tau_H$ usually considered in the literature (see
\cite{cappell-lee-miller}). In our notation, we define it for triples $(L_1,L_2,L_3)$ of Lagrangian subspaces of a
symplectic Hilbert space with almost complex structure $J$, such that $(J L_{i}, L_{j})$ is a Fredholm pair
for all $i,j = 1,2,3$. We set
$\tau_\mu
  (L,L,L):=0$ for some Lagrangian subspace $L$ and use \cite[Formula
  (6.21)]{kirk-lesch} to define $\tau_\mu$ for other triples: If $L_{1,t}$, $L_{2,t}$ and $L_{3,t}$,
  $t\in[0,1]$ are paths of Lagrangian subspaces, such that $(J L_{i,t}, L_{j,t})$ is a Fredholm pair
for all $i,j = 1,2,3$, $t\in[0,1]$, then the (twisted) Maslov triple
  index $\tau_\mu$ is determined by requiring that
\bea
\lefteqn{\tau_\mu(L_{1,1},L_{2,1},L_{3,1}) -
  \tau_\mu(L_{1,0},L_{2,0},L_{3,0})}\\
&=&
\Mas(J L_1,L_2)+\Mas(J L_2, L_3) - \Mas(J L_1, L_3).
\eea
The indices $\tau_\mu$ and $\tau_H$ share some properties. For example $\tau_\mu$ is
additive under direct sums (symplectic additivity). Furthermore the following
properties are an elementary consequence of the above characterization, and can also be found in \cite[Proposition 6.11]{kirk-lesch}.
\begin{lem}\label{tauswitch} Let $L_i$, $i=1,\ldots, 4$ be pairwise Fredholm
  Lagrangians in a Hilbertspace $H$. Then
\bi
\item $\tau_\mu(L_1,L_1,L_2) = \tau_\mu(L_1,L_2,L_2) =
  0$, and
\item $\tau_\mu(L_1,L_2,L_1) = \dim(L_1\cap L_2)$.
\item $\tau_\mu(L_1,L_2,L_3) = \dim(L_2\cap L_3) - \tau_\mu(L_1,L_3,L_2)$.
\ei
\end{lem}

\section{A family of Atiyah--Patodi--Singer boundary
conditions}\label{torusboundaryconditions}

This section introduces a specific family of Atiyah--Patodi--Singer
boundary conditions for the odd signature operator on a manifold with
torus boundary.

The results will be used for the analysis of the spectral flow for
the splitting $M= X \cup_T S$ in Section \ref{SFonClosedMfld},
where $S$ is the solid torus. In Section, \ref{SFonSolidTorus} we
will explicitly compute spectral flow on $S$ with the boundary
conditions developed in this section.

\subsection{A family of flat connections on the torus}
\label{repconn}

Let $\chi_T = \hom(\pi_1 T,SU(2))/ \text{conj}$ be the set of conjugacy classes of
$SU(2)$ representations of the fundamental group $\pi_1T$ of $T$.

The holonomy map gives a
homeomorphism from the gauge equivalence classes of the flat $SU(2)$ connections $\cF_T/\cG_T$ on $T$ to $\chi_T$ \cite[Proposition 2.2.3]{donaldson-kronheimer}. If $A = - i \alpha\, d
m - i \beta\, dl$ with $(\alpha,\beta)\in\R^2$, then
$\text{hol}(A)= \rho_{(\alpha,\beta)}$ is given in
quaternionic notation by
\begin{equation}\label{reponS}
\ba{rcl}
\rho_{(\alpha,\beta)}\co \pi_1 (T) & \to & SU(2)\\
\mu &\mapsto & e^{2\pi i \alpha} \\
\lambda & \mapsto & e^{2\pi i \beta}.\\
\ea
\end{equation}
Notice, that $\R^2\to \chi_T, (\alpha,\beta)\to
[\rho_{(\alpha,\beta)}]$ is a branched cover of $\chi_T$, with
branch points the half integer lattice, which map to central
representations, and covering transformations $(\alpha,\beta)\to
(\pm \alpha + m, \pm \beta + n)$, $(m,n)\in \Z^2$. Each
$(\alpha,\beta)$ then also corresponds to an $SU(2)$ connection
$-i \alpha\, d m - i \beta\, d l$ on $T$. Thus,
 We have a smooth family $\{-i
  \alpha\, d m- i\beta\, d l\}$ of flat connections with holonomy $\rho_{(\alpha,\beta)}$ as in (\ref{reponS}) parametrized by $\R^2$,
which  projects onto $\cF_T/\cG_T = \chi_T$.

\begin{defn}\label{R2family}
Let $a_{\alpha,\beta} := - i \alpha \, dm - i \beta \, dl$. We substitute an index $a_{\alpha,\beta}$ by $(\alpha,\beta)$, for example  $\cH^{0}_{\alpha,\beta}(T;su(2)) =
\cH^{0}_{a_{\alpha,\beta}}(T;su(2))$ and $S_{\alpha,\beta} =
S_{a_{\alpha,\beta}}$.
\end{defn}

The  harmonic forms of $\Delta_{\alpha,\beta} =
\Delta_{a_{\alpha,\beta}}$ are equal to the kernel of
$S_{\alpha,\beta}$ by Lemma \ref{eigcorr}. By $e^{i(2 \alpha m
+2\beta
  l)} \C j$ we denote the $\C$--linear combinations of the function $e^{i(2 \alpha m
+2\beta
  l)}j$. The following
straightforward computation of the harmonic forms on the torus is
left to the reader.

\begin{prop}\ \label{CohoTorus}
\bea
\cH^0_{\alpha,\beta}(T;su(2)) & = &
\begin{cases}
\R i & \text{for }(\alpha,\beta)\notin(\hz)^2,\\
\R i\oplus e^{i(2 \alpha m +2\beta
  l)} \C j & \text{for }(\alpha,\beta)\in (\hz)^2
\end{cases}\\
\cH^1_{\alpha,\beta}(T;su(2)) & = &
\begin{cases}
\R i\, dm \oplus \R i \, dl & \text{for
}(\alpha,\beta)\notin(\hz)^2,\\
\left\{{\setlength{\arraycolsep}{0cm}\begin{array}{c} \R i\, dm \oplus e^{i(2 \alpha m +2\beta
  l)} \C j \, dm\\
{}\oplus \R i\, dl \oplus e^{i(2 \alpha m +2\beta
  l)} \C j \, dl
\end{array}}\right\}
 & \text{for
}(\alpha,\beta)\in (\hz)^2
\end{cases}
\\
\cH^2_{\alpha,\beta}(T;su(2)) & \cong &
\cH^0_{\alpha,\beta}(T;su(2)) \text{ via the Hodge star.}
\eea
\end{prop}

Also observe that $S_{\alpha,\beta}$ preserves the $\R i$ and $\C j$--part of $L^2(\Om^{0+1+2}(T;su(2)))$, and
we can write
$$\ P^{\pm}_{\R i} = P^{\pm}_{\alpha,\beta}
\cap L^2(\Om^{0+1+2}(T;\R i))\text{ and } P^{\pm}_{\alpha,\beta,\C j}  =  P^{\pm}_{\alpha,\beta} \cap L^2(\Om^{0+1+2}(T;\C j)).$$
as well as
$$\cH^{0+1+2}(T;\R i)=\cH^{0+1+2}_{\alpha,\beta}(T;\R i).$$

\subsection{Boundary conditions}\label{boundaryconditions}

We want a path of self-adjoint operators which is continuous in
the graph topology. In order for the
odd signature operator $D_{A_t}$  on $S$ or $X$ to vary continuously and be
(unbounded) self-adjoint, we need a path of boundary conditions which is
continuous in the gap topology. See \cite{booss-lesch-phillips}
for details.

We might like to pick a continuous family of Lagrangians in
$\cH^{0+1+2}_{\alpha,\beta}(T;su(2))$ and extend it to a family of
APS boundary conditions in $\LOT{su(2)}$ by $P^+_{\alpha,\beta}$
or $P^-_{\alpha,\beta}$. Unfortunately, not only does the
dimension of the kernel of the tangential operator jump up at
$(\hz)^2$ (see Proposition \ref{CohoTorus}), but for a smooth path
$\varrho_t$, $t\in [-1,1]$, through $\varrho_0\in (\hz)^2$ we have
$$\lim_{t\to 0^+} P^\pm_{\varrho_t} \neq \lim_{t\to 0^-}
P^\pm_{\varrho_t}.$$ In fact the $\cH^{0+1+2}_{\alpha,\beta}(T;\C
j)$ part of the limits turn out to be orthogonal to each other.

\begin{figure}[ht!]\anchor{r2tilde}\small
\begin{center}
\leavevmode
\includegraphics[scale=.7]{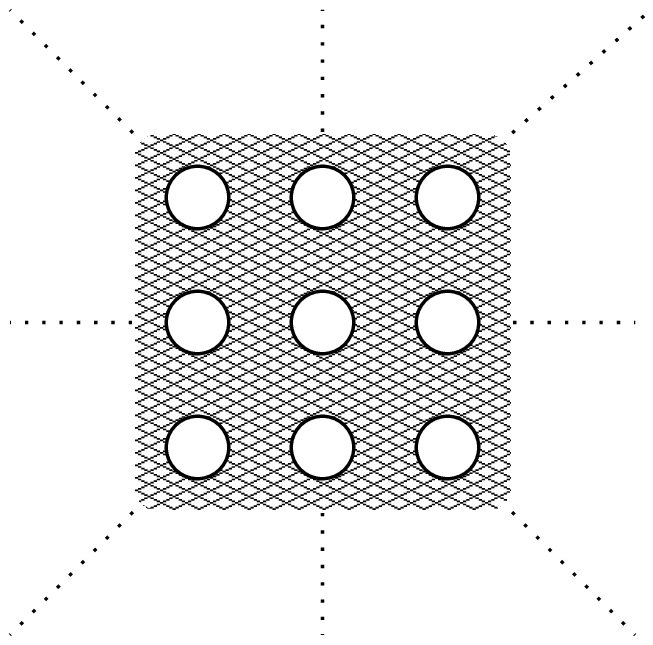}
\end{center}
\caption{$\dover\R^2 \homeo \Rt^2$\label{r2tilde}}
\end{figure}

Thus we introduce the space $\dover\R^2$ shown in 
\figref{r2tilde} to parametrize the Lagrangians. It is $\R^2$ with
open disks of radius $\frac{1}{8}$ removed around all half integer
lattice points with the induced topology. Some people would call
it the real blow-up of the plane at the half integer lattice
points. We will see in Theorem \ref{Ktheta}, why this is a good
parameter space. The advantage of this space is that we can easily
homotop paths of connections together with their boundary
conditions, thus getting a homotopy of paths of self-adjoint
operators. However, it will be more convenient to parametrize
$\dover\R^2$ by the following space.

\begin{defn} \label{homeomorphism}  Let $\widetilde{\R}^2:=\R^2 \times S^1 /\sim$, where $(\alpha,\beta,\theta)
  \sim (\alpha,\beta,1)$ if $(\alpha,\beta)\notin
  (\hz)^2$. We will simply write $(\alpha,\beta,\theta)\in \Rt^2$. Alternatively it is convenient to think of elements in
  $\widetilde\R^2$ as being of the form
$$
\begin{cases}
(\alpha,\beta,\theta)& \text{ if } (\alpha,\beta)\in (\hz)^2
\text{
    and } \theta \in S^1,\\
(\alpha,\beta)& \text{ if } (\alpha,\beta)\notin (\hz)^2.
\end{cases}$$
Denote by $\pi\co \Rt^2 \to \R^2$ the projection
$(\alpha,\beta,\theta)\to (\alpha,\beta)$. Let us define a
bijection $h$ between $\Rt^2$ and $\dover\R^2$.

We will describe what this bijection looks like around the origin.
At all the other half integer lattice points we get a similar
bijection via translation. Away from disks of radius $\frac{1}{4}$
around each half integer lattice point the bijection is the
identity map. Identify $\R^2$ with $\C$ in the usual way. Let
$D\subset\C$ be the disk of radius $\frac{1}{4}$, $\tilde D :=
D\times S^1/\sim$ and $A\subset \C$ the disk of radius
$\frac{1}{4}$ with an open disk of radius $\frac{1}{8}$ around the
origin removed, that is, an annulus. Let $\eta\co \R \to [0,1]$ be a
smooth (cutoff) function with
$$
\eta(t) =
\begin{cases}
0 & t \le \frac{1}{8}\\
\text{a homeomorphism onto }[0,1] & \frac{1}{8} \le t \le \frac{1}{4}\\
1 & t \ge \frac{1}{4}
\end{cases}
$$
Then the bijection $h\co A \to \tilde D$ is given by
$h(z):=(\eta(|z|)\cdot z,\frac{z}{|z|}) \in \tilde D$. We give
$\Rt^2$ the topology that makes the bijection $h\co\dover\R^2 \to
\Rt^2$ into a homeomorphism. Notice that $h$ is a diffeomorphism
away from $(\hz)^2$.
\end{defn}

\begin{thm}\label{Ktheta} Let $(\alpha,\beta)\in \R^2$ and $\theta \in
 S^1 \subset \C$. Then
\be
\item $P^\pm_{\alpha,\beta}\co \R^2-(\hz)^2 \to
  \{\text{closed subspaces of }\LOT{su(2)}\}$ is continuous.
\item Moreover $\lim_{t\to
  0^+} P^\pm_{(\alpha,\beta)\pm t(\re\, \theta,\im\, \theta)}$ exists
  and we can define $$K^\pm_{(\alpha,\beta,\theta)}:=\lim_{t\to
  0^+} P^\pm_{(\alpha,\beta)\pm t(\re\, \theta,\im\, \theta)}/
  P^\pm_{(\alpha,\beta)} \subset \cH_{(\alpha,\beta)}^{0+1+2}(T;\C
  j).$$ Note that $K^\pm_{(\alpha,\beta,\theta)} = 0$ for
  $(\alpha,\beta)\notin (\hz)^2$.
\item $P^\pm_{\alpha,\beta} \oplus K^\pm_{(\alpha,\beta,\theta)}\co \Rt^2 \to
  \{\text{closed subspaces of }\LOT{su(2)}\}$ is continuous.
\item If $(\alpha,\beta) = \frac{1}{2} (r,s)\in
(\hz)^2$, then $K^\pm_{(\alpha,\beta,\theta)} = \la \{ \psi_1^\pm
j,\psi_2^\pm j,\psi_1^\pm k, \psi_2^\pm k \} \ra$
$$\psi^\pm_1 = e^{i(rm+sl)}(1\mp(i \im \theta \, dm- i \re
\theta  \, dl))\leqno{\hbox{where}}$$
 $$\psi_2^\pm = e^{i(rm+sl)}( dm \wedge dl \pm( i
\re \theta  \, dm + i \im \theta \, dl)).\leqno{\hbox{and}}$$ \ee
\end{thm}

Notice from the explicit description of $K^{\pm}$ that
$K^{\pm}_{(\alpha,\beta,\theta)} =
K^{\mp}_{(\alpha,\beta,-\theta)}$. Before we prove this theorem,
an eigenspace decomposition will be useful to compute
$K^\pm_{(\alpha,\beta,\theta)}$ explicitly and study the behaviour
of the family $P^\pm_{(\alpha,\beta)} \oplus
K^\pm_{(\alpha,\beta,\theta)}$ around the half integer lattice. It
is a lengthy but straightforward computation, which we will leave to
the reader.

\begin{prop} \label{eigenvectors} Fix $a= a_{\alpha,\beta}$. We have
  the orthogonal decomposition $$\LOT{su(2)} =
  \overline{\bigoplus_{(r,s)\in\Z^2} E_{r,s}}^{L^2} \oplus
  \overline{\bigoplus_{(r,s)\in\Z^2} E'_{r,s}}^{L^2}$$
where forms in $$E_{r,s} = \{f \sin (rm + sl) + g \cos (rm +
sl) \mid f,g \in \cHT{\R
    i} \}$$ are eigenvectors of $\Delta_{a}$ with eigenvalue $r^2 +
 s^2$, and forms in
 $$E'_{r,s} = \{e^{i(rm
+sl)} f \mid f \in \cHT{\C j}\}$$ are eigenvectors  of
$\Delta_{a}$ with eigenvalue $(r-2\alpha)^2+(s-2\beta)^2$.
$\cHT{\R
    i}$ and $\cHT{\C j}$ are the harmonic forms of the untwisted
Laplacian.
\end{prop}

\begin{proof}[Proof of Theorem \ref{Ktheta}] For the continuity of $P_{\alpha,\beta}^\pm$ away from
  $(\hz)^2$ see
\cite{kirk-klassen2}. We show that $\lim_{t\to
  0^+} P^\pm_{(\alpha,\beta)\pm t(\re\, \theta,\im\, \theta)}$ exists
by explicitly computing $K^\pm_{(\alpha,\beta,\theta)}$ for $(\alpha,\beta) = \frac{1}{2}
(r,s) \in (\hz)^2$.

Let $(\alpha_t,\beta_t)=\frac{1}{2} (r,s) + \frac{1}{2}t (\re \theta,
\im \theta)$, $t> 0$. By Proposition \ref{eigenvectors} $\phi = e^{i(r m + s
  l)} j$ is an eigenvector of $\Delta_{(\alpha_t,\beta_t)} = \Delta_{a_{(\alpha_t,\beta_t)}}$ with eigenvalue $t^2 = (r-2\alpha_t)^2 +
(s-2\beta_t)^2$. Then by Lemma \ref{eigcorr}
$$\phi\pm\frac{1}{t}S_{(\alpha_t,\beta_t)}\phi =
e^{i(rm+sl)}\left(1\mp\frac{1}{t}i(r-2\alpha_t)\,dl \pm
\frac{1}{t}i(s-2\beta_t)\,dm\right)j
$$
is a $\pm t$--eigenvector of $S_{(\alpha_t,\beta_t)}$. This yields
$$\lim_{t\to 0}(\phi\pm\frac{1}{t}S_{(\alpha_t,\beta_t)}\phi) = e^{i(r m + s
l)}(1 \mp( i \im \theta\, dm - i \re \theta\, dl)) .$$
For $\phi =
e^{i(r m + s l)} j\,dm \wedge dl$ we similarly get
$$
\lim_{t\to 0}( \phi\pm\frac{1}{t}S_{(\alpha_t,\beta_t)}\phi )=
e^{i(rm+sl)}( dm \wedge dl \pm( i \re \theta \, dm + i \im \theta \, dl)) j.
$$
Repeating the same computation for $\phi = e^{i(r m + s l)} k$ and
$\phi = e^{i(r m + s l)} k\,dm \wedge dl$ yields a total of $8$
linearly independent eigenvectors, which lie in
either $K^+_{(\alpha,\beta,\theta)}$ or
$K^-_{(\alpha,\beta,\theta)}$. By Proposition
\ref{CohoTorus} the $\C j$ part of $\ker
  S_{(\alpha,\beta,\theta)}$ is $8$--dimensional. Thus the
  $L^2$--span of the above eigenvectors make up $K^+_{(\alpha,\beta,\theta)}$ or
$K^-_{(\alpha,\beta,\theta)}$. This completes the computation of
$K^\pm_{(\alpha,\beta,\theta)}$ and shows that $\lim_{t\to
  0^+} P^\pm_{(\alpha,\beta)\pm t(\re\, \theta,\im\, \theta)}$ exists.

We are left with proving the continuity of $P^+_{\alpha,\beta} \oplus
  K^+_{\alpha,\beta,\theta}$  parametrized by $(\Rt)^2$. Away from
  $(\hz)^2$ the family $P^+_{\alpha,\beta} \oplus
  K^+_{\alpha,\beta,\theta}$ is continuous since
  $K^+_{\alpha,\beta,\theta}=0$. To show
  continuity at half integer lattice points, it suffices to show
  that for any continuous path $\tilde\varrho_t=(\tilde\alpha_t,\tilde\beta_t,\theta_t)$ in $\Rt^2$ that limits to
  $\tilde\varrho(0)=(0,0,\theta_0)$ we have $$\lim_{t\to 0} P^+_{\tilde\alpha_t,\tilde\beta_t}\oplus
  K^+_{\tilde\varrho_t} = P^+_{0,0} \oplus K^+_{0,0,\theta}.$$
For all other half integer lattice points the argument is the
same.
By definition $\tilde\varrho_t$ is continuous in a ball of radius
$\frac{1}{8}$
  around $(0,0)$, if $h^{-1}\circ\tilde\varrho_t = \varrho_t =
  (\alpha_t,\beta_t)$ is continuous in $\R^2$, where
  $h\co\dover\R^2 \to \Rt^2$ is the homeomorphism given in Definition \ref{homeomorphism}. We
  consider two cases:
\be
\item If $|\varrho_t| = \frac{1}{8}$ for $t$ small, then
  $\tilde\varrho_t=(0,0,\theta_t)$ and $\theta_t$ continuous for small
  $t$. Elementary triangle equality arguments applied to
  $P^+_{(0,0)+s(\re \theta, \im\theta)}$, $s$ small, show that $K^+_{0,0,\theta_t}$ is
  continuous.
\item Let $|\varrho _t| \neq \frac{1}{8}$ for $t>0$ small.
  We have $(\tilde\alpha_t,\tilde\beta_t) = (\alpha_t,\beta_t)\eta(\sqrt{\alpha_t^2 + \beta_t^2})$ and
  $\theta_0 = 8(\alpha_0+i\beta_0) \in S^1$. To check that $K^+_{0,0,\theta_0} = \lim_{t\to 0}
  (P^+_{\tilde\alpha_t,\tilde\beta_t}\slash P^+_{0,0})$, observe that
  by Lemma
  \ref{eigcorr} the
  $((2\tilde\alpha_t)^2+(2\tilde\beta_t)^2)$--eigenvector $\phi= j$ of
  $\Delta_{\tilde\alpha_t,\tilde\beta_t}$ yields a $2\sqrt{\tilde\alpha_t^2+\tilde\beta_t^2}$--eigenvector of $S_{\tilde\alpha_t,\tilde\beta_t}$:
\bea \lefteqn{\phi +
\frac{1}{2\sqrt{\tilde\alpha_t^2+\tilde\beta_t^2}}S_{\tilde\alpha_t,\tilde\beta_t}
  \phi
 =  \phi + \frac{2 i \tilde \alpha_t\, dl - 2i \tilde \beta_t\,
  dm}{2\sqrt{\tilde\alpha_t^2 + \tilde \beta_t^2}}}\\
& = & j + \frac{i\alpha_t\eta(\sqrt{\alpha_t^2+\beta_t^2})dl -
i\beta_t \eta(\sqrt{\alpha_t^2+\beta_t^2}) dm}{
  \sqrt{\alpha_t^2\eta(\sqrt{\alpha_t^2 + \beta_t^2})^2+\beta_t^2\eta(\sqrt{\alpha_t^2 + \beta_t^2})^2}}j \\
& = & j + \frac{i\alpha_t dl - i\beta_t dm}{
  \sqrt{\alpha_t^2+\beta_t^2}}j \\
& \stackrel{t\to 0}{\longrightarrow} & j + 8(i\alpha_0 dl -
i\beta_0 dm) j = j - ( i \im
  \theta\, dm   - i \re \theta\, dl)j
\eea The same computation for the other
  $2\sqrt{\tilde\alpha_t^2+\tilde\beta_t^2}$--eigenvectors of
  $S_{\tilde\alpha_t,\tilde\beta_t}$ yields the other basis elements of $K^+_{0,0,\theta}$ as
  given in the statement of Theorem \ref{Ktheta}.
\ee
Thus $P^+_{\alpha,\beta} \oplus
  K^+_{\alpha,\beta,\theta}$ is a continuous family parametrized by $(\Rt)^2$.
\end{proof}

\begin{defn}\label{LambdaL}
For a continuous family of Lagrangians $L$ of $\cH^{0+1+2}(T;\R
i)$ parametrized by a subset $U\subset \Rt^2$, define a family
$\APS^\pm(L)$ of subspaces of $\LOT{su(2)}$ para\-met\-rized by
$U$ as follows
$$
\APS_{(\alpha,\beta,\theta)}^\pm(L) :=
(P^\pm_{\R i} \oplus L_{(\alpha,\beta,\theta)}) \oplus (P^\pm_{\C
j,
  (\alpha,\beta)} \oplus K^\pm_{(\alpha,\beta,\theta)})
$$
\end{defn}

For our application, the family of Lagrangians
$L_{\alpha,\beta,\theta}$ will be independent of
$(\alpha,\beta,\theta)$. Notice that $K_{(\alpha,\beta,\theta)}$
vanishes away from $(\hz)^2$, while we have ``blown up'' the
points of $(\hz)^2$ and removed the singularities that paths
through $(\hz)^2$ encounter. Here is a corollary of Theorem
\ref{Ktheta}, which is important for the following sections.

\begin{thm}\label{LambdaLcontinuous}
$\APS^\pm(L)$ is a continuous family of Lagrangians parametrized
by $\widetilde\R^2$.
\end{thm}

\section{Spectral flow on a closed
$3$--manifold}\label{SFonClosedMfld}

In this section we develop a splitting formula for spectral flow, which expresses
spectral flow of the odd signature operator between flat
connections on a closed $3$--manifold $M$ in terms of spectral flow on
the solid torus $S$ and its complement $X$ with the Atiyah--Patodi--Singer
boundary conditions from Section \ref{torusboundaryconditions}.
Even though Nicolaescu \cite{nico} and Daniel \cite{daniel}
provide us by way of Theorem \ref{nicoextension} with an elegant
expression of spectral flow in terms of the Maslov index of the
respective Cauchy data spaces, it is not immediately applicable to
spectral flow computations, since the Cauchy data spaces
themselves are complicated objects. The purpose of the splitting
formula in Theorem \ref{splittingformula}, however, is to make
computations of spectral flow easier by shifting the problem to
two more tractable ones. In fact, an explicit way to compute
the spectral flow on $S$ is given in Theorem
\ref{SFSolidTorus}. One only needs to do some work computing
the spectral flow on $X$. The main application to keep in mind is
the spectral flow computation of the twisted odd signature
operator between flat connections, whenever it is possible to find
a path between them which is flat on $X$, because then we can use topology to try computing
the spectral flow on $X$.

\subsection{Objective}\label{setup}

The setup in this section is as follows. Let $M = X \cup_T  S$ be
a closed $3$--manifold with $S$ being the solid torus and $T$ the
torus as in Section \ref{conventions}. Let $A_t$ be a path of $SU(2)$ connections on $M$ with the
following properties:
\begin{enumerate}
\item $A_t$ is in cylindrical form and flat in a collar of $T$.
\item \label{rhocondition} $A_t$ restricts to the path $a_{\varrho(t)}$ (see Definition \ref{R2family})
 on $T$ for
  some path $\tilde\varrho$ in $\Rt^2$ with $\pi\circ \tilde\varrho = \varrho$, where $\pi\co \Rt^2 \to
\R^2$ is the projection onto the $\R^2$--factor.
\item $A_0$ and $A_1$ are flat on $M$.
\end{enumerate}

The goal in this section is to find a splitting
formula expressing $\SF(D_{A_t})$ in terms of spectral flow on $X$
and $S$. Notice that, while the spectral flow on
$X$ and $S$ depends on the lift $\tilde\rho$
of $\rho$ in property (\ref{rhocondition}), $A_t$ and
$\SF(D_{A_t})$ are independent of it.

The above assumptions do not limit the applicability of the splitting
formula in Theorem \ref{splittingformula}. Indeed, the spectral flow of the
odd signature operator along a path of $SU(2)$ connections $A_t$ with
flat endpoints depends only on $A_0$ and $A_1$, which we can gauge
transform by some $g_\epsilon$, $\epsilon= 0,1$ so that $g_\epsilon \cdot
A_\epsilon|_{N(T)} = a_{\alpha_\epsilon,\beta_\epsilon}$. The change
in spectral flow is given by Proposition \ref{degreeshift}. We can extend
the path $(1-t) a_{\alpha_0,\beta_0} + t a_{\alpha_1,\beta_1}$ by obstruction theory to a
path of connections on $M$ with endpoints $g_0 \cdot
A_0|_{N(T)}$ and $g_1 \cdot
A_1|_{N(T)}$, since $[0,1]\times
\Omega^*(M;su(2))$ is contractible.

\subsection{The scattering Lagrangian of $D_{A}$ on $S$}

Before we start with the discussion of the spectral flow on
$M$, we analyze the scattering Lagrangian of $D_{A}$ at a flat connection $A$ on
$S$ which restricts to $a_{\alpha,\beta}=-i\alpha\, dm - i
\beta\, dl$ on $T$, because its explicit description plays a central role in the splitting theorem.

The scattering Lagrangian of $D_A$ for $A = -i\beta\,dl$ can be computed
  directly.

Consider a flat $SU(2)$ connection $A$ on $S$
with
  $A|_T = a_{\alpha,\beta}$. One observes that $p\co\pi_1(T) \to \pi_1(S) = \la \pi_1(T) | \mu = 1
  \ra$ is a surjection and $\hol(A) \circ p = \hol(a_{\alpha,\beta})$. Thus we have
  $\hol(a_{\alpha,\beta})(\mu) = 1$, that is, $\alpha \in
  \Z$, and in view of Proposition \ref{scattidentif}
the scattering Lagrangian depends
  only on $a_{\alpha,\beta}$. Therefore we can use a gauge transformation $g$ with $g|_T =
  e^{i\alpha m}$ on the boundary to compute $\cL_{S,g\cdot A} = \ad_g
  \cL_{S,A}$. This yields the following.

\begin{lem}\label{ScatLagOnS} Suppose $A$ is a flat $SU(2)$ connection on $S$
with
  $A|_T = a_{\alpha,\beta}$. Then $\alpha \in \Z$ and
\bi
\item  $\cL_{S,A} = \R i \oplus \R i\, dl$ for $\beta \in \R-\hz$,
\item $\cL_{S,A} = \hat
  \cL_{S,A} \oplus \check \cL_{S,A}$ for $\beta \in \hz$, where
\begin{eqnarray*}
\hat \cL_{S,A} & = & \R i \oplus \R i \, dl\\
\check \cL_{S,A} & = & e^{i(2 \alpha m + 2\beta l)} (\C j \oplus
\C j \, dl)
\end{eqnarray*}
\ei
\end{lem}

Thus we make the following definition.

\begin{defn}\label{hatL}
For $(\alpha,\beta) \in \R^2$ define $\hat \cL_{S,A_{(\alpha,\beta)}}
:= \R i \oplus \R i\, dl$. Note that $\hat \cL_{S,A_{(\alpha,\beta)}}
= \cH^0(T;\R i) \oplus \R
i \, dl \subset \cH^{0+1+2}(T;\R i)$.
\end{defn}

\subsection{A splitting formula for spectral flow}

We will derive the following splitting formula.

\begin{thm} \label{splittingformula} Let $M = X \cup_T  S$ be a closed $3$--manifold with $S$ being the
solid torus and $T$ the torus as in Section
\ref{conventions}. Let $A_t$ be a path of $SU(2)$ connections on
$M$ with the following properties:
\begin{enumerate}
\item $A_t$ is in cylindrical form and flat in a collar of $T$.
\item $A_t$ restricts to the path $a_{\varrho(t)}$ on $T$ for
  some path $\tilde\varrho$ in $\Rt^2$ with $\pi\circ \tilde\varrho = \varrho$, where $\pi\co \Rt^2 \to
\R^2$ is the projection onto the $\R^2$--factor.
\item $A_0$ and $A_1$ are flat on $M$.
\end{enumerate}
Then we have the splitting formula: \bea \SF(D_{A_t}) &\! \! \!  =
&\! \! \!
\SF(D_{A_t}|_S;\APS^+_{\tilde\varrho(t)}(J \hat \cL_S))) +
\SF(D_{A_{t}}|_X;\APS^-_{\tilde\varrho(t)}(\hat \cL_S)))\\
&&\hspace{-0.35cm}{}+\tau_\mu(J\cL_{S,0},K^+_{\varrho_0,i} \oplus J \hat
\cL_{S,0},\cL_{X,0}) -\tau_\mu(J\cL_{S,1},K^+_{\varrho_1,i} \oplus J \hat
\cL_{S,1},\cL_{X,1}).
\eea
\end{thm}

The proof of Theorem \ref{splittingformula} is deferred to the next
section. However, the idea is the following. First we relate the
spectral flow to the Maslov index of the Cauchy data spaces using
Nicolaescu's splitting theorem. Next we use homotopy invariance of
the Maslov index and Nicolaescu's adiabatic limit theorem to find
homotopic paths of Lagrangians, parts of which correspond to a
Maslov index of Cauchy data space with APS boundary conditions.
Now we can use additivity of pairs of paths of Lagrangians and
apply the relative version of Nicolaescu's splitting theorem to
two of the parts, which yields the first two summands of the
formula. The remaining summands will either vanish or
simplify to the two Maslov triple indices involving
$K^\pm_{\alpha,\beta,\theta}$ from Theorem \ref{Ktheta}, the
scattering Lagrangians $\cL_{S,\epsilon} = \cL_{S,A_\epsilon}$ and
$\cL_{X,\epsilon} = \cL_{X,A_\epsilon}$ at $\epsilon = 0,1$, as well as
$\hat{\cL}_{S,\epsilon} = \hat{\cL}_{S,A_\epsilon}$ from Lemma
\ref{ScatLagOnS}.

The situation is particularly nice and straightforward when
\be
\item $D_{A_0}$ and $D_{A_1}$ restricted to $X$ have
non-resonance level $0$, and
\item $\varrho(t) \notin (\hz)^2$ for all $t\in [0,1]$, in which case
  $K^\pm_{(\varrho(t),\theta)}=0$.
\ee
We describe the derivation of our splitting formula without these
  assumptions.

Fix the path $R_t:=\frac{t}{1-t}$ which corresponds
to the parameter that stretches the collar from $0$ to $\infty$.
 Stretching the collar of $S$ yields the adiabatic limit
 $\Lambda^{R_1}_{S,\epsilon} = P^-_{S,\epsilon} \oplus
 \cL_{S,\epsilon}$, because $D_{A_\epsilon}|_S$ has non-resonance level
 $0$ when $\epsilon=0,1$. Theorem \ref{adiabaticlimit} describes $\Lambda^{R_1}_{X,\epsilon}$.

Since $K^\pm_{\alpha,\beta,\theta} \neq
K^\pm_{\alpha,\beta,-\theta}$ for $(\alpha,\beta)\in (\hz)^2$ (see Theorem \ref{Ktheta}), we will use the APS boundary conditions
$\APS_{\tilde\varrho(t)}^-(\hat \cL_S)$ and $\APS_{\tilde\varrho(t)}^+(J \cL_S)$ for our splitting formula with the notation introduced
in Definition \ref{LambdaL} and Lemma \ref{ScatLagOnS}, which were
shown to be continuous for a continuous lift $\tilde\varrho$ of the
path $\varrho$ from $\R^2$ to $\Rt^2$.

We will need to ``rotate''
the adiabatic limit of $D_{A_0}|_S$ and $D_{A_1}|_X$ to our preferred
APS boundary conditions.

\begin{defn}\label{KWpaths}  We introduce the paths $P^+_{1,\nu} \oplus L_{W,t} \oplus
d_{a_1}(E^-_{1,\nu}) \oplus L_{X,t}$ for $X$ and $P^-_0 \oplus L_{S,t}$
for $S$ as follows:
\begin{enumerate}
\item
Let $L_{X,t}$ be a path of Lagrangians in $\HT{su(2)}$ from the
scattering Lagrangian $\cL_{X,1}$ to $J \hat \cL_{S,1} \oplus
K^+_{\varrho_1,i}$.
\item $L_{W,t} = e^{J \pi t}W_{a_1} \oplus J(W_{a_1}^\perp)$ is a
  path of Lagrangians in $d_{a_t}(E^+_{a_1,\nu})\oplus
  d^*(E^-_{a_1,\nu})$ from $W_{a_1} \oplus
J(W_{a_1}^\perp)$ to $d^*(E^-_{a_1,\nu}) = J(W_\alpha) \oplus
J(W_a^\perp) \subset P^+_1$.
\item For $\epsilon=0,1$ let $L_{S,\epsilon,t} =
  \hat\cL_{S,\epsilon} \oplus \check L_{S,\epsilon,t}$
  be a path of Lagrangians in
  $\HT{su(2)}$, where $\check L_{S,\epsilon,t}$ is an arbitrary path from
$\check \cL_{S,\epsilon}$ to $K^-_{\varrho_\epsilon,i}$. As in Lemma
\ref{ScatLagOnS}, $\check
  \cL_{S,\epsilon}$ denotes the $\C j$ part of
  $\cL_{S,\epsilon}$. Set $L_{S,t} = L_{S,0,t}$.
\end{enumerate}
\end{defn}

Then it is straightforward to check that the composition of paths
$(\cM_i(t),\cN_i(t))$ given in \tabref{pathsext} is homotopic to
the path of pairs of Cauchy data spaces
$(\Lambda_{S,t},\Lambda_{X,t})$. The table shows the $i$--th paths
$\cM_i$ and $\cN_i$ in the second and fifth column. The third and
fourth column give the endpoints of $\cM_i$ and $\cN_i$ as a
reference.

\newcommand{\rb}[1]{\raisebox{1.2ex}[-1.2ex]{$\displaystyle{#1}$}}
\begin{table}[ht!]\anchor{pathsext}\small
\begin{center}
\leavevmode
$
\setlength{\arraycolsep}{0.05cm}
\begin{array}{c|c|c|c|c}
i\vrule width 0pt depth 5pt\relax & \text{paths } {\cM}_{i}(t) &
\multicolumn{2}{c|}{\text{Endpoints
of }{\cM}_{i}\text{ and }{\cN}_{i}} & \text{paths } {\cN}_{i}(t)\\
\hline\hline
 & & & & \\
\cline{1-2}\cline{5-5} & & \rb{\Lambda_{S,0}} & \rb{\Lambda_{X,0}} & \\
\cline{3-4} \rb{1} & \rb{\Lambda_{S,0}^{R_t}} & & & \rb{\Lambda^{R_t}_{X,0}} \\
\cline{1-2}\cline{5-5} & & \rb{P^-_0\oplus \cL_{S,0}} &  \rb{\Lambda^\infty_{X,0}} & \\
\cline{3-4} \rb{2} & \rb{P^-_0\oplus L_{S,t}} & & & \rb{\Lambda^{R_{1-t}}_{X,0}} \\
\cline{1-2}\cline{5-5} & & \rb{\APS^-_{\tilde\varrho(0)}(\hat \cL_{S})} & \rb{\Lambda_{X,0}} & \\
\cline{3-4} \rb{3} & \rb{\APS^-_{\tilde\varrho(t)}(\hat \cL_S)} & & &
\rb{\Lambda_{X,t}} \\
\cline{1-2}\cline{5-5} & & \rb{\APS^-_{\tilde\varrho(1)}(\hat \cL_{S})} & \rb{\Lambda_{X,1}} & \\
\cline{3-4} \rb{4} & \rb{\text{constant}} & & & \rb{\Lambda_{X,1}^{R_t}}\\
\cline{1-2}\cline{5-5} & & \rb{\APS^-_{\tilde\varrho(1)}(\hat \cL_S)} & \rb{\Lambda^\infty_{X,1}} & \\
\cline{3-4} \rb{5} &  \rb{\text{constant}} & & & \rb{P^+_{a_1,\nu}
  \oplus L_{W,t} \oplus d_{a_1}(E^-_{a_1,\nu}) \oplus L_{X,t}} \\
\cline{1-2}\cline{5-5} & & \rb{\APS^-_{\tilde\varrho(1)}(\hat \cL_S)}
&
\rb{\APS_{\tilde\varrho(1)}^+(J \hat \cL_{S})} & \\
\cline{3-4} \rb{6} &  \rb{\APS^-_{\tilde\varrho(1-t)}(\hat \cL_S)} &
& &
\rb{\APS_{\tilde\varrho(1-t)}^+(J \hat \cL_{S})}\\
\cline{1-2}\cline{5-5} & & \rb{\APS^-_{\tilde\varrho(0)}(\hat \cL_S)}
&
\rb{\APS^+_{\tilde\varrho(0)}(J \hat \cL_{S})} & \\
\cline{3-4} \rb{7} &  \rb{P^-_0\oplus L_{S,1-t}} & & & \rb{\text{constant}}\\
\cline{1-2}\cline{5-5} & & \rb{P^-_0\oplus \cL_{S,0}} &
\rb{\APS^+_{\tilde\varrho(0)}(J \hat \cL_{S})} & \\
\cline{3-4} \rb{8} & \rb{\Lambda_{S,0}^{R_{1-t}}} & & & \rb{\text{constant}}\\
\cline{1-2}\cline{5-5} & & \rb{\Lambda_{S,0}} & \rb{\APS^+_{\tilde\varrho(0)}(J \hat \cL_{S})} & \\
\cline{3-4} \rb{9} & \rb{\Lambda_{S,t}} & & &\rb{\APS^+_{\tilde\varrho(t)}(J \hat \cL_{S})}\\
\cline{1-2}\cline{5-5} & & \rb{\Lambda_{S,1}} &
\rb{\APS_{\tilde\varrho(1)}^+(
J \hat \cL_{S})} & \\
\cline{3-4} \rb{10} & \rb{\Lambda^{R_t}_{S,1}} & & &
\rb{P^+_{a_1,\nu}
  \oplus L_{W,1-t} \oplus d_{a_1}(E^-_{a_1,\nu}) \oplus L_{X,1-t}}\\

\cline{1-2}\cline{5-5} & & \rb{\Lambda^{\infty}_{S,1}} & \rb{\Lambda^\infty_{X,1}} & \\
\cline{3-4} \rb{11} & \rb{\Lambda^{R_{1-t}}_{S,1}} & & &  \rb{\Lambda_{X,1}^{R_{1-t}}} \\
\cline{1-2}\cline{5-5} & & \rb{\Lambda_{S,1}} &  \rb{\Lambda_{X,1}} & \\
\end{array}
$\\[3ex]
\end{center}
\caption{\label{pathsext}The paths homotopic to $\Lambda_{S,t}$
and $\Lambda_{X,t}$
  broken up into pieces}
\end{table}
If $\varrho_\epsilon\notin(\hz)^2$ for $\epsilon = 0,1$, then we have $K^\pm_{\varrho_\epsilon, i} = 0$, thus
  $$\tau_\mu(J\cL_{S,0},K^+_{\varrho_0,i} \oplus J \hat
\cL_{S,0},\cL_{X,0}) = \tau_\mu(J\cL_{S,1},K^+_{\varrho_1,i} \oplus J \hat
\cL_{S,1},\cL_{X,1}) =
  0.$$ This yields the following.

\begin{cor}\label{hzsplittingformula} We make the same assumptions as in Theorem
  \ref{splittingformula}. Then,
 if $\varrho_\epsilon\notin(\hz)^2$ for $\epsilon = 0,1$ we get:
$$
\SF(D_{A_t}) = \SF(D_{A_{t}}|_S;\APS^+_{\tilde\varrho(t)}(J \hat \cL_S)) +
\SF(D_{A_{t}}|_X;\APS^-_{\tilde\varrho(t)}(\hat \cL_S)).
$$
\end{cor}

A way to compute $\SF(D_{A_t}|_S;\APS^+_{\tilde\varrho(t)}(J \hat \cL_S))$
will be given in Theorem \ref{SFSolidTorus}. We will compute
$\SF(D_{A_{t}}|_X;\APS^-_{\tilde\varrho(t)}(\hat \cL_S))$ for our main application in
Section \ref{SFonTorusBundle}, where $M$ is a torus-bundle over
$S^1$ and $A_t$ is flat on $X$.

\subsection{Proof of Theorem \ref{splittingformula}}

By Lemma \ref{intersectionstretching} the Maslov indices $\Mas(\cM_i,\cN_i)$ vanish for
$i=1,4,8,11$. By Theorem \ref{nicorelative} we have $$\Mas(\cM_9,\cN_9) =
\SF(D_{A_{t}}|_S;\APS^+_{\tilde\varrho(t)} (J \hat \cL_{S})).$$
Since
$\APS^-(\hat \cL_{S,t}) = (\APS^+_{\tilde\varrho(t)}(J \hat
\cL_{S}))^\perp$ for all $t$, the Maslov index $\Mas(\cM_6,\cN_6)$ vanishes.
Again by Theorem \ref{nicorelative} we have $$\Mas(\cM_3,\cN_3) =
 \SF(D_{A_t}|_X; \APS^-(\hat \cL_S)).$$
Let us focus our attention on the rest of the paths, namely paths
2, 5, 7 and 10.

We may homotop $\cM_2$ to the composition of the two
  paths
\bea
\cM_{2a} & = & P^-_0\oplus L_{S,0,t}\\
\cM_{2b} & = & P^-_0\oplus K^+_{\varrho_0,i}\oplus \hat \cL_{0}. \eea
and $\cN_2$ to the composition of the two
  paths
\bea \cN_{2a} & = & P^+_{0,\nu} \oplus (W_0 \oplus J(W^\perp_0))
\oplus
d^*(E^-_{0,\nu}) \oplus \cL_{X,0}\\
\cN_{2b} & = & \Lambda_{X,0}^{R_{1-t}}. \eea By Lemma \ref{intersectionstretching} $\Mas(\cM_{2b},\cN_{2b}) = 0$. Then by our
choice of $L_{S,0,t}$ we get
$$\Mas(\cM_2,\cN_2) =\Mas(\cM_{2a},\cN_{2a}) = \Mas (L_{S,0,t},\cL_{X,0}).$$
We can directly check:
$$\Mas(\cM_7,\cN_7) = -\Mas (L_{S,0,t}, K^+_{\varrho_0,i} \oplus J \hat \cL_{S,0})).$$
Then Lemma \ref{tauswitch} and linearity under the
 direct sum of triples yield
 \begin{align*}
\lefteqn{\Mas(\cM_2,\cM_2) + \Mas(\cM_7,\cN_7)}\\
& = \tau_\mu(J(K^-_{\varrho_0,i} \oplus \hat
\cL_{S,0}),\cL_{X,0},K^+_{\varrho_0,i} \oplus J \hat
\cL_{S,0}) - \tau_\mu(J\cL_{S,0},\cL_{X,0},K^+_{\varrho_0,i} \oplus J \hat \cL_{S,0})\\
& =  \dim (\cL_{X,0} \cap (K^+_{\varrho_0,i} \oplus J \hat \cL_{S,0})) - \tau_\mu(J\cL_{S,0},\cL_{X,0},K^+_{\varrho_0,i} \oplus J \hat \cL_{S,0})
\\
& =  \tau_\mu(J\cL_{S,0},K^+_{\varrho_0,i} \oplus J \hat \cL_{S,0}, \cL_{X,0})
\end{align*}
We get a similar expression for the paths 5 and 10. The path $\cM_{10}$ can be homotoped to the
  composition of the two paths
\bea
\cM_{10a} & = & \Lambda^t_{S,1}\\
\cM_{10b} & = & \Lambda^\infty_{S,1} =  P^-_0\oplus \cL_{S,1} \eea
and $\cN_{10}$ to the composition of the two
  paths
\bea
\cN_{10a} & = & \APS^+_{\tilde\varrho(1)}(J \hat \cL_{S})\\
\cN_{10b} & = & P^+_\nu \oplus L_{W,t} \eea  Again by Lemma
\ref{intersectionstretching} we get $\Mas(\cM_{10a},\cN_{10a}) =
  0$. Thus
\bea \Mas(\cM_{10},\cN_{10}) & = & \Mas (P^-_1\oplus \cL_{S,1},
P^+_{1,\nu}
  \oplus L_{W,1-t} \oplus d_{a_1}(E^-_{1,\nu}) \oplus L_{X,1-t})\\
& = & \Mas (E^-_1\oplus \cL_{S,1}, L_{W,1-t} \oplus d_{a_1}(E^-_{1,\nu})
  \oplus L_{X,1-t})\\
& = & \Mas (E^-_1, L_{W,1-t} \oplus d_{a_1}(E^-_{1,\nu})) + \Mas
(\cL_{S,1}, L_{X,1-t}) \eea while \bea \Mas(\cM_{5},\cN_{5}) & = &
\Mas (P^-_1\oplus K^-_{1,i} \oplus \hat \cL_{S,1}, P^+_{1,\nu}
  \oplus L_{W,t} \oplus d_{a_1}(E^-_{1,\nu}) \oplus L_{X,t})\\
&  = &  \Mas (E^-_1, L_{W,t} \oplus d_{a_1}(E^-_{1,\nu})) + \Mas
(K^-_{1,i}
  \oplus \hat \cL_{S,1}, L_{X,t})
\eea Thus
$$ \Mas(\cM_{5},\cN_{5}) + \Mas(\cM_{10},\cN_{10}) = \Mas (\cL_{S,1},
L_{X,1-t}) + \Mas (K^-_{1,i}  \oplus \hat \cL_{S,1},
  L_{X,t}).
$$
Notice that since $L_{S,1,t}\circ L_{S,1,1-t}$ and $L_{X,t}\circ
L_{X,1-t}$ are both homotopic to constant paths, we can compute
\bea\lefteqn{\Mas(\cM_{5},\cN_{5}) + \Mas(\cM_{10},\cN_{10})}\\
 & = & \Mas (\cL_{S,1},
L_{X,1-t}) + \Mas (K^-_{1,i} \oplus \hat \cL_{S,1},
   L_{X,t})\\
& = & - \Mas (L_{S,1,t},\cL_{X,1}) + \Mas (L_{S,1,t},
   K^+_{\varrho_1,i}\oplus J\hat \cL_{S,1}).
\eea
Just like for the paths 2 and 7 this simplifies to
\begin{align*}
\lefteqn{\Mas(\cM_5,\cM_5) + \Mas(\cM_{10},\cN_{10})}\\
& =  \tau_\mu(J(K^-_{\varrho_1,i} \oplus \hat
\cL_{S,1}),K^+_{\varrho_1,i} \oplus J \hat \cL_{S,1},\cL_{X,1}) - \tau_\mu(J\cL_{S,1},K^+_{\varrho_1,i} \oplus J \hat
\cL_{S,1},\cL_{X,1}) \\
& =  - \tau_\mu(J\cL_{S,1},K^+_{\varrho_1,i} \oplus J \hat
\cL_{S,1},\cL_{X,1}).
\end{align*}
This shows that $\sum_{i=1}^{9}
\Mas(\cM_i,\cN_i)$ simplifies to the desired formula and completes the proof.\qed

\section{Spectral flow on the solid torus}\label{SFonSolidTorus}

In our main application discussed in Section \ref{SFonTorusBundle} we
are interested in computing spectral flow of the twisted odd
signature operator along a path of $SU(2)$
connections on torus bundles $M$ over $S^1$. The main
tool is our splitting
formula for spectral flow from Theorem \ref{splittingformula}, by which we ultimately need to compute spectral flow
on a solid torus $S \subset M$ and on its complement $M - S$ with
certain APS boundary conditions.

While the spectral flow computation on $M - S$ certainly depends on $M$ itself,
we will discuss spectral flow on $S$
separately, as it might be interesting and applicable in other
settings.

\subsection{Objective}

Let $S = D^2 \times S^1$ be the solid torus and $T = \partial S$
its boundary. Let $A_{\alpha,\beta}$ be the family of
connections on $S$ as described in Definition \ref{connectionsOnS} and let
$A_t = A_{\alpha_t,\beta_t}$ be a path of connections on the solid torus $S$, so that $A_0$ and $A_1$ are flat connections.
Boundary conditions are necessary to make the odd signature
$D_{A_t}$ on $S$ (unbounded) self-adjoint.
The Atiyah--Patodi--Singer (APS) boundary
conditions $\APS^+_{(\alpha,\beta,\theta)}(L)$ discussed in Section \ref{torusboundaryconditions} are suitable for our purpose.
The goal is to compute the spectral flow $\SF(D_{A_{\varrho_t}};\APS^+_{\tilde\varrho_t}(L))$ of
$D_{A_t}$ on $S$, where $\varrho = \pi \circ \tilde\varrho$.

Since $A_{\alpha,\beta}$ with $\alpha \in \Z$ covers all gauge
equivalence classes of flat connections on $S$, the restriction to a
path in $A_{\alpha,\beta}$ is suitable for its use in a
splitting formula, as we discussed in Section \ref{setup}.

\subsection{A smooth family of connections on $S$ parametrized by
  $\R^2$}\label{FamilyOnS}

We would like to extend the family of flat connections on $T$
given in Definition \ref{R2family} to a family $A_{(\alpha,\beta)}$
(smoothly parametrized by $\R^2$) of connections on $S$, so that
\be
\item $A_{(\alpha,\beta)}$ is flat for
$\alpha \in \Z$,
\item the
restriction of $A_{(\alpha,\beta)}$ to a collar
$N(T)$ of
$T$ in $S$ is $a_{(\alpha,\beta)} := - i\alpha\, dm - i
\beta\, dl \in \cA_{N(T)}$.
\ee

Notice that $i\beta\, dl$ makes sense globally as a (flat) connection on $S$ for all $\beta \in \R$, whereas $i\alpha\,
dm$ does not for $\alpha \neq 0$. Thus we can try to construct the
family, so that
\be
\item[(3)] $A_{(\alpha,\beta)}$ is equal to $-i \beta\, dl$ for
$\alpha=0$ and $\beta\in \R$.
\ee

Since $\pi_2 SU(2) = 0$, we can find a gauge transformation
$\fa$ on $S$, such that in the collar $N(T)$ of $T$
\begin{eqnarray*}
\fa|_{N(T)}\co \quad N(T) & \to & SU(2)\\
(n,e^{i m},e^{i  l}) & \mapsto & e^{i m}.
\end{eqnarray*}
The gauge
  transformation $\fa$ can easily be constructed factoring through
  $D^2$. It will be convenient to fix $\fa(ne^{im},e^{il}) = q(n)
  e^{im}+ \sqrt{1-(q(n))^2}j$ for a smooth non-decreasing
  cutoff function $q\co [0,1] \to [0,1]$ with $q(n)=0$ for $n$ near $0$ and $q(n)=1$.
Since
\begin{eqnarray*}
\fa|_{N(T)} \cdot a_{(\alpha,\beta)} & =& \fa a_{(\alpha,\beta)} \fa^{-1} +
\fa d(\fa^{-1})
=   a_{(\alpha,\beta)} + e^{im}
\tfrac{\partial}{\partial m} e^{-i m} dm\\
& = & a_{(\alpha,\beta)} - i dm = a_{(\alpha+1,\beta)},\nonumber
\end{eqnarray*}
we get a smooth family of connections
$A_{(\alpha,\beta)} := -\fa^{\alpha} \cdot i\beta\, dl$ on $S$, $(\alpha,\beta)\in
\Z \times \R$, with the desired property
$$
A_{(\alpha,\beta)}|_{N(T)} = \fa^{\alpha}|_{N(T)} \cdot a_{(0,\beta)} = a_{(\alpha,\beta)}.
$$
This can be easily extended to a smooth family of connections
parametrized by $\R^2$ with
$A_{(\alpha,\beta)}|_{N(T)} = a_{(\alpha,\beta)}$. One possible way is
the following:

\begin{defn}\label{connectionsOnS}
Let $\eta\co S \to [0,1]$ be a smooth cutoff function with $\eta(x)=1$ for
$x\in N(T)$ and $\eta(x)=0$ near the core of $S$, which may as well factor through $D^2$. Let
$\{\tau_r\co\R \to \R\}_{r\in\Z}$ be a partition of unity
subordinate to $\{(r-1,r+1)\}_{r\in\Z}$. For $(\alpha,\beta)\in\R^2$ define the family
$$A_{(\alpha,\beta)} := \sum_{r\in\Z} \tau_r(\alpha)\, \fa^{r}\cdot
(\eta\, i(r-\alpha)\, dm - i \beta\, dl).$$
\end{defn}

Notice that $\eta\, i(r-\alpha)\, dm
- i \beta\, dl$ is an $SU(2)$ connection  on $S$ (in normal form) for
$(\alpha,\beta)\in \R^2$ and $r\in\Z$. Thus
$A_{(\alpha,\beta)}$ is a smooth family of $SU(2)$ connections on $S$.
It is easy to check, that this family has the following
properties:
\be
\item $A_{(\alpha,\beta)} = \fa \cdot A_{(\alpha+1,\beta)}$.
\item
$
A_{(\alpha,\beta)}|_{N(T)} = a_{(\alpha,\beta)}.
$
\item $A_{(\alpha,\beta)} = - \fa^{-\alpha}\cdot i\beta\, dl$ for $\alpha\in\Z$, which
  is flat.
\ee

\subsection{Computation of spectral flow on the solid torus}

In this section all connections and odd signature operators are
considered on the solid torus only. Everything that follows depends on
a continuous family of
Lagrangians $L_{\alpha,\beta,\theta}$ in $\cH^{0+1+2}(T;\R
 i)$ parametrized by $\Rt^2$.

To avoid technical difficulties we are going to assume that $L_{\alpha,\beta,\theta}$ is  a family of Lagrangians in $\cH^{0+1+2}(T;\R
 i)$, which is transverse to $\hat\cL_{S,A_{(\alpha,\beta)}}$ (see Definition \ref{hatL}) for all
 $(\alpha,\beta,\theta) \in \Z\times \R \times \{ \pm i\} \subset
 \Rt^2$. Later we will fix a specific family of Lagrangians $L_{\alpha,\beta,\theta}$ in view of
 the splitting formula in Theorem \ref{splittingformula}.

\begin{defn} Define $\SF(\tilde\varrho):=\SF(D_{A_{\pi\circ\tilde\varrho(t)}},\APS^+_{\tilde\varrho(t)}(L))$
 as a function of paths $\tilde\varrho$ in $\Rt^2$, where $\pi\co\Rt^2 \to
 \R^2$ is the projection mentioned in Definition \ref{homeomorphism}.
\end{defn}

Observe that $\SF$ is additive under compositions of paths in $\Rt^2$
and that $\SF$ is a homotopy invariant rel endpoints.

\begin{figure}[ht!]\anchor{cycle}\small
\begin{center}
\leavevmode
\psfrag{2}{$2$}
\psfrag{-4}{$-4$}
\psfrag{4}{$4$}
\psfrag{-12}{$-12$}
\psfrag{12}{$12$}
\psfrag{ZxRx{i}}{$\Z\times\R\times \{\pm i\}$}
\psfraga <-2pt,2pt> {hz}{$(\hz)^2\times S^1$}
\includegraphics[scale=.7]{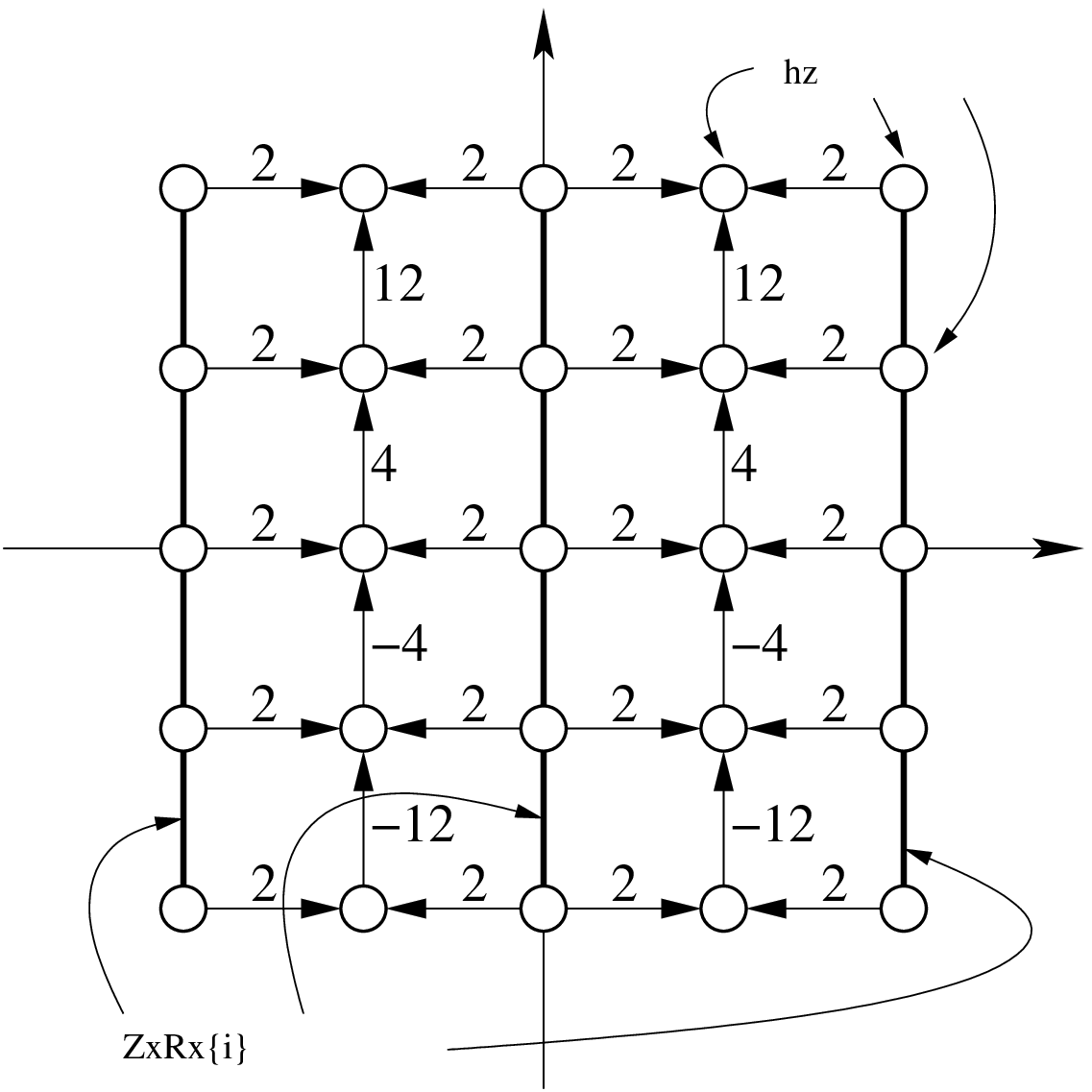}
\end{center}
\caption{\label{cycle}Cycle in $\Rt^2$}
\end{figure}

\begin{lem}\label{SFequalintersection}
There is a cohomology class $z\in H^1(\Rt^2,\Z\times\R\times \{ \pm i \})$, such that for a path $\tilde\varrho$ in
$\Rt^2$ that starts and ends in $\Z\times \R\times \{\pm i\} \subset
\Rt^2$ the (mod 0) spectral flow equals $z(u)$, where
$u:=[\tilde\varrho([0,1])]\in H_1(\Rt^2,\Z\times\R\times \{ \pm i \})$. Note, that $\Z\times
 \R \times \{ \pm i \}$ corresponds to the thickened vertical lines in
 \figref{cycle}.
\end{lem}

\begin{proof}
By defining $\zeta(\sum_i a_i \sigma_i):=\sum_i a_i\SF(\sigma_i)$ for
singular $1$--simplices $\sigma_i$ and integers $a_i$, we get a map
$\zeta\co S_1(\Rt^2) \to \Z$.

Now suppose $(\alpha,\beta)\in\tilde\varrho([0,1]) \subset \Z \times \R \times \{
 \pm i\}$. Then $A_{\tilde\varrho}$ is
  a path of flat connections on $S$. We have
  $$\ker(H^1(S;su(2)_{\hol(A_{\alpha,\beta})})\to
  H^1(T;su(2)_{\hol(A_{\alpha,\beta})}))
  =0.$$ Then Proposition
  \ref{sescor} and  $L_{\alpha,\beta} \cap \cL_{S,A_{\alpha,\beta}}
  = 0$ imply that $D_{A_{\alpha,\beta}}$ (with boundary conditions ${\APS^+_{\alpha,\beta}(L)}$)
  has no kernel. Furthermore by comparing
  Theorem
  \ref{Ktheta} and Lemma \ref{ScatLagOnS} we see that
  $K^+_{\alpha,\beta,\pm i} \cap \check \cL_{S,A_{\alpha,\beta}} = 0$
  when $(\alpha,\beta)\in \Z \times \frac{1}{2} \Z$.
This implies that $\SF(\tilde\varrho) = 0$.

Thus $\zeta$ descends to a map $\bar\zeta\co S_1(\Rt^2)/S_1(\Z\times \R \times \{\pm
i\}) \to \Z$, that is $\bar\zeta \in S^1(\Rt^2,\Z\times \R \times \{\pm
i\})$.

It is straightforward to show using the homotopy invariance and the
additivity of $\SF$ that $\bar\zeta$ is a cocycle. We also observe
that by definition a coboundary in $S^1(\Rt^2,\Z\times \R \times \{\pm i\})$ vanishes on paths $\tilde\varrho$ in $\Rt^2$ with endpoints in $\Z \times
\R \times \{ \pm i\}$. Thus we have found the desired cohomology
class $z = [\bar\zeta]$.
\end{proof}

Equivalently we may say (by Poincar\'e duality) that there exists an infinite
  homology class in $H_1(\Rt^2-(\Z \times \R\times \{ \pm i\})))$, such that the
  intersection number with (the homology class representing) the path is
  equal to the spectral flow.

Fix  $L_{\alpha,\beta,\theta} := J
\hat\cL_{S,A_{(\alpha,\beta)}} =  \R i \, dm\oplus \R i\, dm \wedge
dl = \R i\, dm \oplus \cH^2(T;\R i)$. Our goal is to compute the
(infinite) cycle corresponding to the cohomology class $z$ in Lemma
\ref{SFequalintersection}, therefore extending \cite[Theorem 6.4]{kirk-klassen}.

\begin{thm}\label{SFSolidTorus}
If $\tilde\varrho$ starts and ends in $\Z
  \times \R \times \{ \pm i \}$, then
  $\SF(\tilde\varrho)$ is given by the intersection number $\tilde\varrho
  \cdot z$ where $z$ is the cycle shown in 
  \figref{cycle}.
\end{thm}

\begin{figure}[ht!]\anchor{positiveintersection}\small
\begin{center}
\leavevmode
\psfrag{a}{$a$}
\psfrag{rho0}{$\tilde\varrho(0)$}
\psfrag{rho1}{$\tilde\varrho(1)$}
\includegraphics[scale=.7]{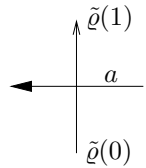}
\end{center}
\caption{\label{positiveintersection} An example with intersection number $a$}
\end{figure}

The sign convention for computing the intersection number is determined by
the example in \figref{positiveintersection}.

The lengthy proof of Theorem \ref{SFSolidTorus} is deferred to the
next section. However, here is a short outline of the
proof. Consider a path of connections $A_t$ between flat connections on a certain closed
manifold $M=S \cup_T X$, which restricts
to a path $A_{\varrho(t)}$ in our family of connections on $S$ defined
in Definition
\ref{connectionsOnS} and $\varrho(t)$ lifts to a path
$\tilde\varrho(t)$ in $\Rt^2$.
An example of such a path $\tilde\varrho(t)$ is shown in
\figref{SF4}. Now consider a particular gauge transformation $g$ on
$M$. Suppose that $\varrho' = \pi \circ \tilde\varrho'$ and
$(g \cdot A_t) |T = A_{\varrho'(t)} | T$.
 Then we are going to
compare the spectral flow of $g\cdot A_t$ with the spectral flow of $B_t$ given by
$$
B_t =
\begin{cases}
g\cdot A_t & \text{when restricted to } $X$\\
A_{\varrho'(t)} & \text{when restricted to } $S$,
\end{cases}
$$
using some properties of a particular family of gauge
transformations as well as Proposition
\ref{degreeshift}.

By repeating the above process a few times and applying the splitting
 formula for spectral flow in Theorem \ref{splittingformula} for different paths $\varrho$
 we will learn how the chains of the desired (infinite) cycle relate to each other, because the terms
 of the spectral flow on $X$ cancel. This yields the
 cycle with three unknowns $a$,$b$ and $c$. The
 coefficients in \figref{cycle} can then be computed using some lens space
examples and by applying the splitting formula for spectral flow in
Theorem \ref{splittingformula} once more.

Before we embark on a proof of Theorem \ref{SFSolidTorus} we adapt
 some of the definitions and results in
 \cite{boden-herald-kirk-klassen} to our needs.

Consider the following group of gauge transformations:
$$
\cG_{\text{nf}} = \{\text{smooth } g\co S \to SU(2)\approx S^3 \mid
g|_{N(T)} (ne^{im},e^{il}) = e^{i\alpha m + i\beta l}, \alpha,\beta \in \Z\}.
$$
We have $H_3(S^3,S^1;\Z) = \Z$ by the long exact sequence of the
pair $(S^3,S^1)$, and $H_3(X,T;\Z) = \Z$ for
$\partial X = T$. Thus there is a well-defined
degree for maps $(X,T) \to (S^3,S^1)$, particularly for elements of $\cG_{\text{nf}}$.

\begin{lem}[Lemma 4.1, \cite{boden-herald-kirk-klassen}] \label{whenhomotopic} Let $g,g'\in \cG_{\text{nf}}$. Then
  $g$ and $g'$ are homotopic if and only if $g|_T = g'|_T$ and $\deg(g)=\deg(h)$.
\end{lem}

Consider the following gauge transformations on $S =
  D^2 \times S^1 \to SU(2)$:
\be
\item $\fa(ne^{im},e^{il}) = q(n) e^{im}+ \sqrt{1-(q(n))^2}j$ from
  Section \ref{FamilyOnS},
\item $\fb(ne^{im},e^{il}) = e^{il}$,
\item $\fc(ne^{im},e^{il})$ is a gauge transformation of degree $1$
  with $\fc |_{N(T)} \equiv 1$.
\ee
The exact description of these gauge transformations is not relevant for us. However, we will need to exploit some of their
properties.
\begin{thm}[Lemma 4.3, Theorem 4.4, \cite{boden-herald-kirk-klassen}]\label{degreetheorem}
We have $[\fa,\fb] = \fc^{-2}$
up to homotopy and
$\deg(\fa^a\fb^b\fc^c) = c-ab.$
\end{thm}

We are going to use one additional gauge transformation on $S$. We set $\fd
\equiv j$. Notice, that even though $\fd \notin \cG_{\text{nf}}$, we have $\fd g \fd^{-1}
\in \cG_{\text{nf}}$ for all $g\in \cG_{\text{nf}}$. Furthermore the
following observation is noteworthy because we will need it later.
\begin{lem} \label{extradegree}$\deg(\fa \fd \fa
\fd^{-1}) = 0$.
\end{lem}

\begin{proof}
We have $H_3(S^3-p,S^1) = 0$ for $p \notin S^1$ by the long exact
 sequence of the pair $(S^3-p,S^1)$. Thus, if a map $g\co  (X,T) \to
(S^3,S^1)$ with $\partial X = T$ misses a point $p\notin S^1$, then it
 is a composition of maps
$$
(X,T) \to (S^3-p,S^1)\hookrightarrow
(S^3,S^1),
$$ and therefore $\deg (g) = 0$ by  functoriality of $H_3$.

Since $\fa$ is homotopic to $\fa'(ne^{im},e^{il}) = ne^{im} +
\sqrt{1-n^2} j$, we have for any $(n e^{im},e^{il}) \in
S$ $$\fa' \fd \fa' \fd^{-1}(ne^{im},e^{il}) =
2n^2-1 + 2n\sqrt{1-n^2}e^{im}j \neq \frac{1}{\sqrt{2}} +
\frac{1}{\sqrt{2}}i,$$ which implies $\deg(\fa \fd \fa
\fd^{-1}) = \deg(\fa' \fd \fa'
\fd^{-1}) =0$.
\end{proof}

The following useful fact follows immediately from the relative
Mayer--Vietoris sequence and the long exact sequence of the pairs $(M,T)$ and $(S^3,S^1)$.

\begin{lem}\label{degreesplitting}
Let $M =  S \cup_T X$. If $g\co  M \to S^3$ with $g(T) \subset S^1$, then
$\deg (g) = \deg(g|_S)
  + \deg(g|_X)$.
\end{lem}

We want to close this section with a straightforward Maslov index computation (like the one in the proof of Lemma
  \ref{CircleMas}).

\begin{lem} \label{FlowFlatSolid} Let $(\alpha,\beta) \in \Z \times \hz$. Then for
  small enough $\epsilon>0$ and varying $t\in [-\epsilon, \epsilon]$
 we have $$
\Mas(\hat\cL_{\alpha,\beta},K^+_{\alpha,\beta,\theta e^{\pi i t}}) =
\begin{cases} 2 & \text{if }\theta=\pm 1,\\
0 & \text{otherwise.}
\end{cases}
$$
\end{lem}

\subsection{Proof of Theorem \ref{SFSolidTorus}}

\begin{figure}[ht!]\anchor{lensspace}\small
\begin{center}
\leavevmode
\psfraga <0pt, 2pt> {hz}{$(\hz)^2\times S^1$}
\psfrag{Flat on S'}{Flat on $S'$}
\psfrag{Flat on S}{Flat on $S$}
\includegraphics[scale=.7]{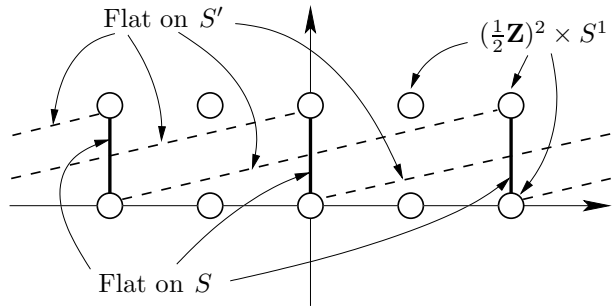}
\end{center}
\caption{Flat connections $S$ and $S'$\label{lensspace}}
\end{figure}

For simplicity consider the lens space $M = L(4,1) = S \cup_h S'$, where $$h\co
\partial S=T \to \partial S', \ (e^{im},e^{il}) \mapsto
(e^{im + il},e^{i4m + i3l}).$$ Then we have a family  flat
connections on $M$ which restricts to noncentral connections on
$T$. \figref{lensspace} shows which flat connections of the
family of connections $a_{\alpha,\beta}$ on $T$ extend to flat
connections on $S$ and $S'$. The intersection of the lines
corresponding to these flat connections on $T$ correspond to flat
connections on $M$.

Since $\pi_1(SU(2)) = \pi_2(SU(2))=0$, there is no obstruction of
extending any $g\co  S \to SU(2)$ to $M$. Denote the
extensions of $\fa$, $\fb$ and $\fc$ by $\bar\fa$, $\bar\fb$ and $\bar\fc$
respectively. We can simply extend $\fd$ to $\bar\fd = j$.

\begin{figure}[ht!]\anchor{SF4}\small
\begin{center}
\leavevmode
\psfrag{0}{$(0,0)$}
\psfrag{1/2}{$(0,\frac{1}{2})$}
\psfrag{1}{$(1,0)$}
\psfrag{a}{$a$}
\psfrag{a'}{$a'$}
\psfrag{beta}{$\fa$}
\psfrag{j}{$\fd$}
\psfrag{rho}{$\tilde\varrho$}
\psfrag{rho'}{$\tilde\varrho'$}
\includegraphics[scale=.7]{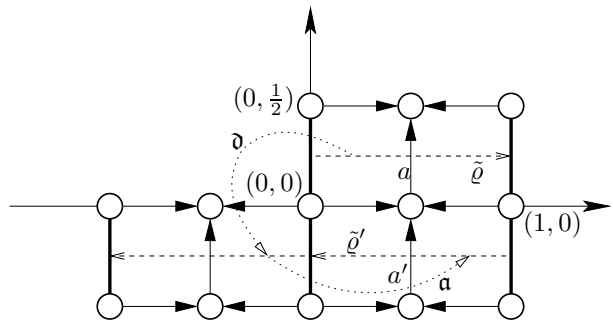}
\end{center}
\caption{Computing $a'=-a$\label{SF4}}
\end{figure}

Let $\tilde\varrho(t)$ be a (straight line) path in $\Rt^2$ which starts at
$(0,\frac{1}{4})\in \Rt^2$ and ends at $(1,\frac{1}{4})\in \Rt^2$ as shown in
\figref{SF4}. Now let $A_t$ be a fixed but arbitrary path of connections on $M$,
which is flat at the endpoints and equals $A_{\varrho(t)}$ for
$\varrho(t) = \tilde \varrho(t)$ when restricted to $S$ as in
Definition \ref{connectionsOnS}. Consider the gauge transformation
$g=\bar \fa \bar
\fd$. Define a path of connections on $M$ by
\begin{equation}\label{Bt}
B_t =
\begin{cases}
g\cdot A_t & \text{when restricted to } $X$\\
A_{\varrho'(t)} & \text{when restricted to } $S$,
\end{cases}
\end{equation}
where $\tilde\varrho'(t)=-\tilde\varrho(t)+(1,0)$ and $\varrho'(t) = \pi\circ\tilde\varrho'(t)$. Notice that $A_{\varrho'(t)}|_T= g\cdot
A_{\varrho(t)}|_T$.

Let $\SF(A_t) := \SF(D_{A_t})$. If $\APS_t$ is a path of APS boundary conditions on $X$ let $\SF_X(A_t;\APS_t):=\SF(D_{A_t}|_X;\APS_t)$.
Occasionally we will also write $\SF(A_0,A_1)$ instead of $\SF(A_t)$
to emphasize that spectral flow on $M$ only depends on the endpoints. Let us compare $\SF(A_t)$
 with $\SF(B_t)$ by utilizing Proposition \ref{degreeshift}. By
 Definition \ref{connectionsOnS} we have $$B_0|_S = \fa \cdot (-
A_0)|_S = \fa \fd \cdot
A_0|_S$$ and $$B_1|_S = \fd \cdot A_0|_S = \fd \fa^{-1} \cdot A_1|_S =
\fd\fa^{-1}\fd^{-1}\fa^{-1}\cdot (\fa\fd \cdot A_1)|_S.$$
Notice that $B_1|_X = \bar\fa \bar\fd \cdot A_1|_X$. Consider the gauge transformation
$$
g'(x) =  \begin{cases} 1 & \text{ if } x\in X\\
\fa\fd\fa\fd^{-1}(x) & \text{ if } x\in S.
\end{cases}
$$
By exploiting the homotopy invariance and the additivity of the
spectral flow,
Proposition \ref{degreeshift} as well as Lemmas
 \ref{degreesplitting} and \ref{extradegree} (in this order) we get
\bea
\SF(A_t)& = & \SF(g\cdot A_t) = \SF(\bar\fa \bar\fd \cdot
A_0,\bar\fa\bar\fd \cdot A_1) =
\SF(B_0, B_1) + \SF(B_1,
\bar\fa \bar\fd \cdot A_1) \\
& =& \SF(B_0, B_1) +
8 \, \deg(g') = \SF(B_0, B_1) +
8 \, \deg(\fa\fd\fa\fd^{-1})\\
& = & \SF(B_t) +
8 \, \deg(\fa\fd\fa\fd^{-1}) =   \SF(B_t).
\eea
Let us now apply the splitting formula in Corollary \ref{hzsplittingformula} to both sides of the equation and
observe that for a path of connections $A_t$, a gauge
  transformation $g$  and a path of APS boundary conditions $\APS_t$ we have
$$
\SF_X(A_t;\APS_t) = \SF_X(g\cdot A_t;\ad_g \APS_t).
$$
In particular $A_t|_T = a_{\tilde\varrho(t)}$ and $g \cdot A_t|_T
= a_{\tilde\varrho'(t)}$ imply
$$
\SF_X(A_t;\APS^-_{\tilde\varrho(t)}(\hat\cL_S)) = \SF_X(g\cdot
A_t;\APS^-_{\tilde\varrho'(t)}(\hat\cL_S)).
$$
Thus the spectral flow terms for $X$ vanish, and $a' = -a$ in \figref{SF4}
because of
$$\SF_S(A_{\varrho(t)};\APS^+_{\tilde\varrho(t)}(J\hat\cL_S)) = \SF_S(A_{\varrho'(t)};\APS^+_{\tilde\varrho'(t)}(J\hat\cL_S)).
$$
\begin{figure}[ht!]\anchor{SF4c}\anchor{SF2a}\small
\begin{center}
\begin{minipage}[t]{6.5cm}
\begin{center}
\leavevmode
\psfrag{0}{$(0,0)$}
\psfrag{1/2}{$(0,\frac{1}{2})$}
\psfrag{1}{$(1,0)$}
\psfrag{a}{$a$}
\psfrag{a'}{$a'$}
\psfrag{alpha}{$\fb$}
\psfrag{rho}{$\tilde\varrho$}
\psfrag{rho'}{$\tilde\varrho'$}
\includegraphics[scale=.7]{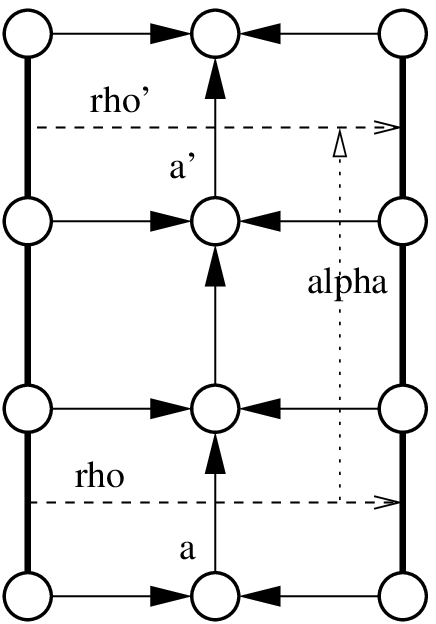}
\caption{Computing $a'=16+a$\label{SF4c}}
\end{center}
\end{minipage}
\begin{minipage}[t]{6cm}
\begin{center}
\leavevmode
\psfrag{b}{$b$}
\psfrag{b'}{$b'$}
\psfrag{c}{$c$}
\psfrag{c'}{$c'$}
\psfrag{alpha}{$\fb$}
\psfrag{rho}{$\tilde\varrho$}
\psfrag{rho'}{$\tilde\varrho'$}
\includegraphics[scale=.7]{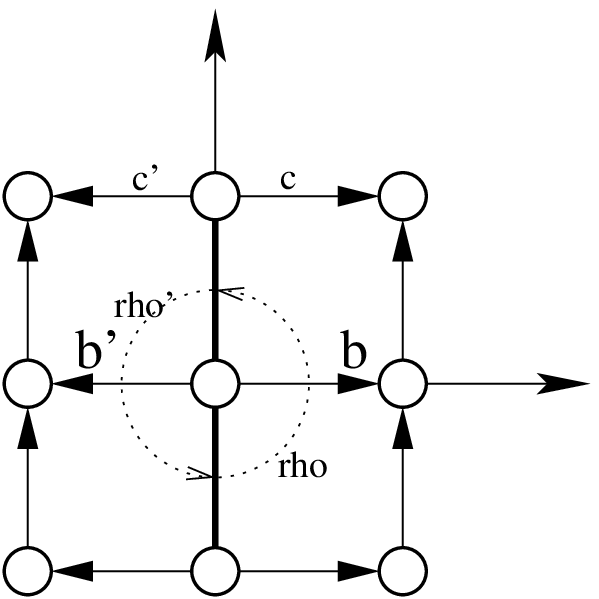}
\caption{Computing $b'=b$\label{SF2a}}
\end{center}
\end{minipage}
\end{center}
\end{figure}
Now let $k,l\in \Z$ and consider the straight line path $\tilde\varrho(t)$ which
starts at $(k,\frac{l}{2}+\frac{1}{4})$ and ends at $(k,\frac{l}{2}+\frac{5}{4})$, as
shown in \figref{SF4c}. As before, let $A_t$ be a fixed but arbitrary path of connections on $M$,
which equals $A_{\varrho(t)}$ for $\varrho(t)=\pi\circ\tilde\varrho(t)$ when restricted to $S$ as in
definition \ref{connectionsOnS}. Consider the gauge transformation
$g=\bar\fb$. Using this data, define a path of connections $B_t$ on $M$ as in
(\ref{Bt}), where $\tilde\varrho'(t)=\tilde\varrho(t)+(0,1)$ and
$\varrho(t) = \pi \circ\tilde \varrho'(t)$. Notice that $A_{\varrho'(t)}|_T= g\cdot
A_{\varrho(t)}|_T$.

We have $$B_0|_S = \fb \cdot A_0|_S$$ and $$B_1|_S = \fa\fb \cdot A_0|_S =
\fa\fb\fa^{-1}\fb^{-1}\cdot (\fb \fa \cdot A_0)|_S =
\fa\fb\fa^{-1}\fb^{-1} \cdot (\fb \cdot A_1)|_S.$$
Consider the gauge transformation
$$
g'(x) =  \begin{cases} 1 & \text{ if } x\in X\\
(\fa\fb\fa^{-1}\fb^{-1})^{-1}(x) & \text{ if } x\in S.
\end{cases}
$$
Then after invoking Theorem \ref{degreetheorem} we get
\bea
\SF(A_t)& = & \SF(g\cdot A_t) = \SF(\bar\fb \cdot A_0,\bar\fb \cdot  A_1) =
\SF(B_0, B_1) + \SF(B_1,
\bar\fb \cdot A_1) \\
& =& \SF(B_0, B_1) + \deg(g')
= \SF(B_0, B_1) +
\deg((\fa\fb\fa^{-1}\fb^{-1})^{-1})\\
& = &\SF(B_t) +
\deg(\fc^2) = \SF(B_t) + 16.
\eea
Application of Corollary \ref{hzsplittingformula} then gives
that $a' = a + 16$ in \figref{SF4c}.

Similarly if for the same path $\varrho$ we consider
$\tilde\varrho'(t)=\tilde\varrho(t)+(1,0)$, $\varrho'(t)=\pi\circ\tilde\varrho'(t)$, the gauge transformation $g=\bar\fa$, and we
define $B_t$ as in (\ref{Bt}), then we have $B_t = g\cdot A_t$ and thus
$\SF(A_t) = \SF(B_t)$.

A similar computation gives us the relationship between coefficients for the horizontal
simplices. Let $\tilde\varrho(t)$ be a path going half around the origin as shown in
\figref{SF2a}, for example $\tilde\varrho(t) = (0,0,\frac{1}{4}e^{\pi i
  (t-\frac{1}{2})})$. Consider $\tilde\varrho'(t) = -
\tilde\varrho(t)$ and $\varrho'(t)=\pi\circ\tilde\varrho'(t)$. Let $g
= \bar\fd$ and  define $B_t$ as in (\ref{Bt}). Since $B_t = \bar\fd
\cdot A_t$, we immediately get $\SF(A_t) = \SF(B_t)$.
The Corollary \ref{hzsplittingformula} yields $b =
b'$ in \figref{SF2a}. Similarly, if $g = \bar\fa$ or $g=\bar\fb$ we have $B_t = g \cdot
A_t$, and we get $c = c'$ in \figref{SF2a}.

\begin{figure}[ht!]\anchor{cycleabc}\small
\begin{center}
\leavevmode
\psfrag{b}{$b$}
\psfrag{a}{$a$}
\psfrag{c}{$c$}
\psfrag{-a+16}{$16-a$}
\psfrag{a-16}{$a-16$}
\psfrag{-a}{$-a$}
\psfrag{ZxRx{i}}{$\Z\times\R\times \{\pm i\}$}
\psfrag{hz}{$(\hz)^2\times S^1$}
\includegraphics[scale=.7]{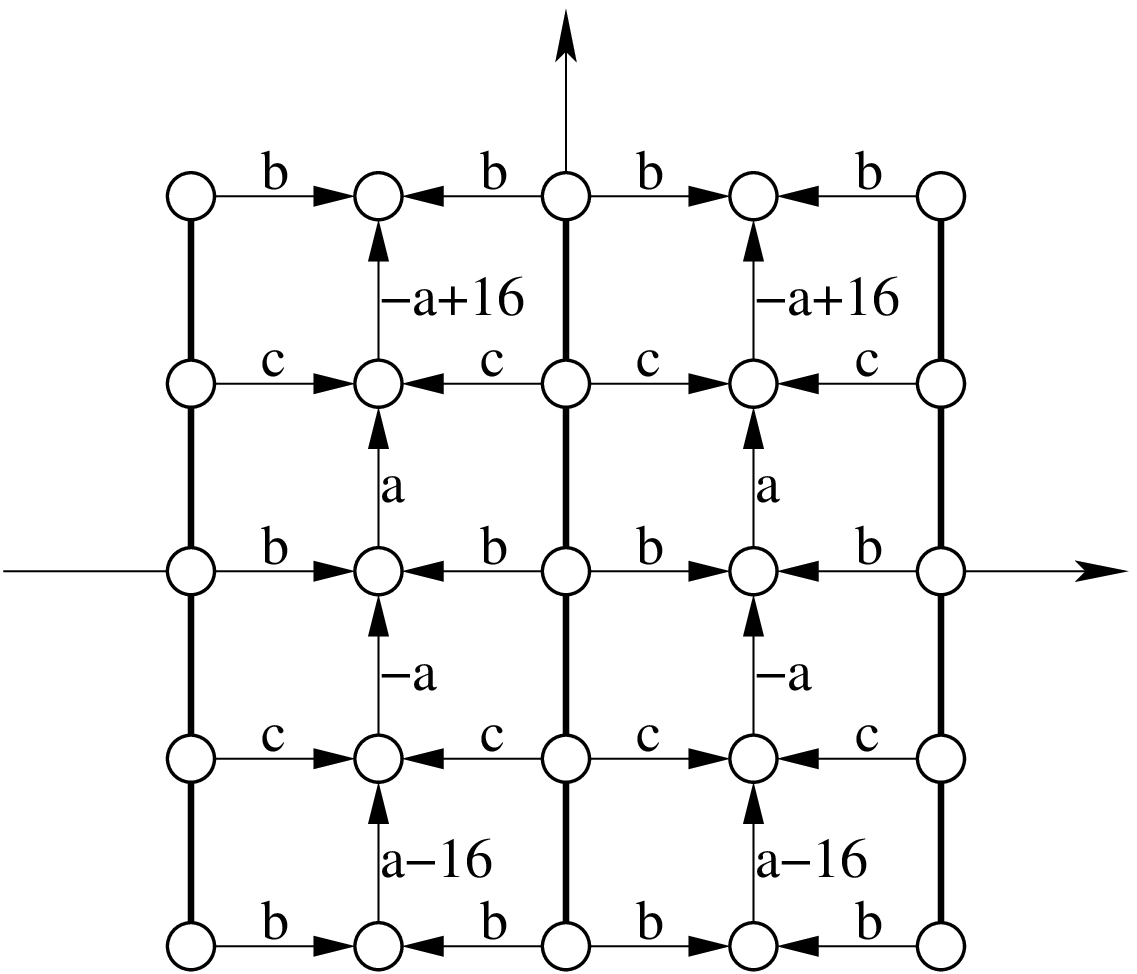}
\end{center}
\caption{\label{cycleabc}Cycle in $\Rt^2$}
\end{figure}

Thus we are left with determining the coefficients $a$, $b$ and $c$ in
our cycle in \figref{cycleabc}.

To find the values of $b$ and $c$ we refer to Lemma
\ref{FlowFlatSolid}. For $a$ we consider the $3$--sphere $M$  constructed
by gluing two solid tori $S$ and
$S'$ with the same orientation along the
boundary using the orientation reversing homeomorphism on the boundary
$(e^{im},e^{il}) \mapsto (e^{il},e^{im})$. For $X=M-S$ we get a similar
statement as in Lemma \ref{FlowFlatSolid}. Then we can employ the splitting formula in
Theorem \ref{splittingformula} for some paths
$\tilde\varrho(t)$ and $\tilde\varrho'(t)$ with
\bi
\item $\tilde\varrho(0)=(0,0)$,
\item $\tilde\varrho(t)=\tilde\varrho'(t)$ for $t\in
[0,\frac{1}{4}]\cup[\frac{3}{4},1]$,
\item $\tilde\varrho(t)=(\frac{1}{2},0,e^{4\pi i t})$ and $\tilde\varrho'(t)$ constant for $t\in
[\frac{1}{4},\frac{3}{4}]$.
\ei
This yields
\bea
0& =& \SF(A_{\tilde\varrho(t)})-\SF(A_{\tilde\varrho'(t)})\\
&=& \SF_S(A_{(\frac{1}{2},0,e^{2\pi i t})};\APS^+(J\hat\cL_S)) +
\SF_{X}(A_{(\frac{1}{2},0,e^{2\pi i t})};\APS^-(\hat\cL_S))\\
&=& (2-a+2-a)+(2+2),
\eea
and shows that
$a$ has to equal
$4$.\qed

\section{Spectral flow on a torus bundle over the circle}\label{SFonTorusBundle}

We are
interested in computing the spectral flow of the odd signature
operator coupled to a path of $SU(2)$ connections on a torus bundle over $S^1$,
because it addresses \cite[Conjecture 5.8]{jeffrey}, the
missing piece in her work on Witten's $3$--manifold invariants
\cite{jeffrey}.

Jeffrey considered a mapping torus $M$ over the torus $T$ and assumed
that its monodromy matrix $B$ has $|\tr B| \neq 2$.
Her conjecture, based on physical reasoning, implies that the spectral
flow of the odd signature operator between
irreducible flat $SU(2)$ connections is
$0 \bmod 4$, which shall be confirmed in this section.

Lisa Jeffrey's conjecture for the case $G=SU(2)$ needs some interpretation
when the trace of the monodromy matrix is zero as well as when one considers
representations of $\pi_1(M)$ whose restriction to $T$ is central (this is equivalent to the
condition in the footnote to
\cite[Conjecture 5.8]{jeffrey} for $G=SU(2)$). Thus we are not able to give a detailed analysis of Lisa
Jeffrey's conjecture in this work. However, it will be the subject of future research.

\subsection{Irreducible $SU(2)$ representations of $\pi_1 M$}

Fix an orientation for $T$. Let $m\co T\nach T$ be an orientation
preserving homeomorphism, $M=T\times I/(m(x),1)\sim
(x,0)$ be its
mapping torus and fix a base point $(*,0)$ in $M$. Consider the isomorphism $m_*\co H_1 (T;\Z) \nach
H_1 (T;\Z)$ induced by $m$ on homology, and fix a meridian $x$ and
longitude $y$ for $T$ so that $H_1 (T;\Z) \cong \pi_1 T = \la x,y, |
[x,y]
\ra \cong \Z^2$ and so that $dx\wedge dy$ is in the orientation class,
when we consider $x,y \in H_1 (T;\Z)\cong H^1 (T;\Z)$. After identifying $x= {1 \choose 0}$ and $y=
{0\choose 1}$ we can write $m_*$ as an unimodular matrix. If $m_*(x)
= ax+cy$ and $m_*(y) = bx + dy$, then the monodromy matrix is
$$
m_*=B=\left[ \ba{cc}
a & b\\
c & d
\ea \right] \in SL_2(\Z).
$$
We will henceforth assume that $|\tr B|\neq 2$. This is equivalent to
$\det(B\pm I)\neq 0$ and implies $c \neq 0$ and $b\neq 0$.

The fundamental group of $M$ is an HNN extension of $\pi_1 T$, where
$\tau$ is the loop from $(*,0)$ to $(m(*),1)$:
$$
\pi_1M=\{x,y,\tau | [x,y], \tau
x\tau^{-1}=x^ay^c,\tau y\tau^{-1}=x^by^d\}.
$$
Recall that we have the identification of $SU(2)$ with the unit quaternions.

Let $\phi = (\phi_1,\phi_2) \in \R^2$. \cite[Proposition 5.5]{kirk-klassen} states that the homomorphism
\begin{equation}\label{rhophi}
\ba{rcl}
\rho_{\phi}\co  \la \{ x,y,\tau\} \ra & \nach & SU(2)\\
\tau &\mapsto & j\\
x & \mapsto & e^{2\pi i\phi_1}\\
y & \mapsto & e^{2\pi i\phi_2}
\ea
\end{equation}
factors through $\pi_1M$ if and only if $\phi(B+I)\in
\Z^2$. Observe that a representation $\rho_\phi$ is reducible if and only if
  $\phi \in (\hz)^2$. Furthermore representations $\rho_\phi$ and $\rho_\psi$ are conjugate if and
only if $\phi = \pm \psi + \theta$ for some $\theta \in \Z^2$. Together with \cite[Corollary 7.2]{kirk-klassen} this implies the
following.
\begin{lem}\label{knot} Any two irreducible
  $SU(2)$ representations $\rho$ and $\rho'$ of $\pi_1(M)$ are conjugate to $\rho_{\phi}$
  and $\rho_{\psi}$ respectively so that
  $\rho_{(1-t)\phi + t\psi}$ as in (\ref{rhophi}) is a path of irreducible $SU(2)$ representations
  of
$$\pi_1(M - N_\gamma) = \la x,y,\tau | [x,y], \tau x^p y^q \tau^{-1} =
x^{pa+qb}y^{pc+qd}\ra,
$$
where $N_\gamma$ is a tubular
  neighborhood of the curve $\gamma=px+qy$ on $T = T \times \{0\}$ and $(p,q)$ is a relatively prime
pair satisfying
\begin{equation}\label{paircondition}
(\phi-\psi)(B+I){p \choose q}=0.
\end{equation}
\end{lem}

Let the {\em longitude} $\lambda$ be the curve on $T$
parallel to $\gamma$ and the {\em meridian} $\mu$ the curve that bounds a disk in
$N_\gamma$. It is not hard to see that that $\mu = x^{ra+sb}y^{rc+sd}\tau
x^{-r}y^{-s}\tau^{-1}$ and $\lambda  =
x^py^q$ are words representing the meridian and longitude. See \cite[Figure 3]{kirk-klassen}. Then for $\phi= (\phi_1,\phi_2)\in\R^2$ we have
\bea
\rho_\phi(\mu) & = & e^{2\pi i \alpha_\phi}\\
\rho_\phi(\lambda) & = & e^{2\pi i \beta_\phi},
\eea
where
$\alpha_\phi := \phi (B+I){r \choose s}$ and $\beta_\phi := \phi {p
  \choose q}$. Thus a representation $\rho_\phi$ of $\pi_1 (M-N_\gamma)$ restricts to a central representation on the boundary if and only if
  $\alpha_\phi\in\hz$ and $\beta_\phi\in\hz$.

\subsection{An example}\label{exple}

It is conjectured by Kirk and Klassen in
\cite{kirk-klassen}, that one can always find some $\phi$ and $\psi$
so that the entire path $\rho_{t\phi + (1-t)\psi}$ is noncentral when
restricted to $\partial (M-N_\gamma)$. The following example shows that
this is false. While Kirk and Klassen's work  does not apply here, the
present paper does apply to compute the spectral flow.
%If this had been true, then their
%work would have covered all cases of paths between irreducible flat
%connections.
% and it is indeed necessary to consider paths that go
%through representations which are central on $\partial
%(M-N_\gamma)$.

Consider $B=\left[ \begin{array}{cc} 5 & 2 \\ 2 & 1 \end{array}
\right]$. Since
$\det(B+I) = 8$, conjugacy classes of $SU(2)$ representations of
$\pi_1(M)$ are uniquely represented by all $\rho_\phi$ with $\phi \in
(\tfrac{1}{8} \Z)^2 \cap ([0,\tfrac{1}{2}] \times [0,\tfrac{1}{2}]
\cup (\tfrac{1}{2},1) \times (0,\tfrac{1}{2}))$ for which $\phi (B+I)
\in \Z$. There are only two conjugacy classes of irreducible
representations. They are represented by $\rho_\phi$ and $\rho_\psi$ for $\phi =
(\tfrac{3}{4},\tfrac{1}{4})$ and $\psi =
(\tfrac{1}{4},\tfrac{1}{4})$.

Any $\rho_{t\phi + (1-t)\psi}$ for $t\in [0,1]$
determines a representation of $\pi_1(M-N_\gamma)$, where $\gamma = px +
qy$ is such that $(p,q)$ is a relatively prime pair satisfying
equation (\ref{paircondition}). In this case we
have $\phi (B+I) = (5,2)$ and $\psi (B+I) = (2,1)$. Thus we
can choose $(p,q)= \pm (1,-3)$. Conjugates $\rho_{\phi'}$ and
$\rho_{\psi'}$ of $\rho_\phi$ and
$\rho_\psi$ give us a different $(p,q)$. For example if we choose $\phi' = \phi$ and
$\psi' = -\psi$, then $(p,q) = \pm(3,-7)$.

\begin{prop}\label{example} Let $B=\left[ \begin{array}{cc} 5 & 2 \\ 2 & 1 \end{array}
\right]$, $\phi =
(\tfrac{3}{4},\tfrac{1}{4})$ and $\psi =
(\tfrac{1}{4},\tfrac{1}{4})$. Let $\phi' = \pm \phi + \eta$ and $\psi' = \pm \psi +
  \theta$ for $\eta,\theta \in \Z^2$ and $\gamma = px + qy$ a knot
  satisfying \ref{paircondition} for $\phi'$ and $\psi'$. Then both $\rho_{\phi'}$ and $\rho_{\psi'}$ are central
  when restricted to $\partial (M-N_\gamma)$.
\end{prop}

\begin{proof} We have $(\phi-\psi)(B+I) = (3,1)$ and $(\phi+\psi)(B+I) =
  (7,3)$. Furthermore for $\theta = (\theta_1,\theta_2)\in \Z^2$
  we get $\theta (B+I) = (6\theta_1+2\theta_2, 2\theta_1+2\theta_2) =
  2(3\theta_1+\theta_2,\theta_1+\theta_2)$. Similarly $\eta(B+I) =
  2(3\eta_1+\eta_2,\eta_1+\eta_2)$. It follows that $(\phi'-\psi') (B+I)$ is a
  pair of odd integers. Since $(p,q)$ are required to be relatively
  prime, one of them must be odd. Since $(\phi'-\psi') (B+I) {p\choose
  q} = 0$, both $p$ and $q$ must be odd. This implies that
  $\beta_{\phi'}$ and $\beta_{\psi'}$ are half
  integers. Since $\alpha_{\phi'}$ and $\alpha_{\psi'}$ are also
  integers,
$\rho_{\phi'}|_{\partial (M-N_\gamma)}$ and $\rho_{\psi'}|_{\partial
  (M-N_\gamma)}$ are central.
\end{proof}

A similar example is given by
$B  = \left( \begin{array}{cc} 3 & 4 \\ 2 & 3 \end{array}
  \right)$ with $\phi= (\frac{1}{4},0)$, $\psi=(\frac{1}{4},\frac{1}{2}).$

In addition to the above examples there are also paths
with representations in the interior which are central on $\partial
(M-N_\gamma)$. The example $$\textstyle B  = \left( \begin{array}{cc} 9 & 4 \\ 2 & 1 \end{array}
  \right) \text{ with } \{\phi= (\frac{1}{6},\frac{1}{6}),\
  \psi=(\frac{2}{3},\frac{1}{6})\} \text{ or } \{\phi= (\frac{5}{6},\frac{1}{3}),\ \psi=(\frac{1}{3},\frac{1}{3})\}$$ are particularly interesting, because all conjugate
choices of $\rho_\phi$ and $\rho_\psi$ seem to still have
representations in the interior which are central on the boundary,
though no proof has been found.

\subsection{Computation of the spectral flow}

We compute the mod 4 spectral flow of the  odd signature operator
on $M$ coupled to a path of $SU(2)$ connections $A_t$, where $A_0$
and $A_1$ are flat and irreducible. We want to apply the splitting formula in
Theorem \ref{splittingformula}. By Lemma \ref{knot} we may assume
that $\rho_{\phi_0} = \hol(A_0)$ and
$\rho_{\phi_1} = \hol(A_1)$, and that $\rho_{\phi_t}$ with $\phi_t := (1-t)\phi_0 + t
  \phi_1$ for $t\in [0,1]$ is a path of irreducible representations of $\pi_1(M-
N_\gamma)$ for some curve $\gamma$. The holonomy maps flat
connections on $M-N_\gamma$ to representations of $\pi_1(M-N_\gamma)$. By \cite[Theorem 4.1]{fine-kirk-klassen} there is local splitting of this holonomy
map. Thus
$\rho_{\phi_t}$ lifts to a path of irreducible connections on $X:=M- N_\gamma$
with holonomy $\rho_{\phi_t}$, which extends to a path $A_t$ of
connections on all of $M$ via the
family of connections $A_{\alpha,\beta}$ on $S := N_\gamma$ given in Definition
\ref{connectionsOnS}. Notice that $A_t|_S = A_{\alpha_t,\beta_t}$
where $(\alpha_t,\beta_t) := (\alpha_{\phi_t},\beta_{\phi_t})$. The
slope of $(\alpha_t,\beta_t)$ corresponds to
$\theta^\star :=
\frac{(\alpha_1-\alpha_0)+(\beta_1-\beta_0)i}{\sqrt{(\alpha_1-\alpha_0)^2+(\beta_1-\beta_0)^2}}\in
S^1$.

Kirk and Klassen computed the spectral flow of $D_{A_t}$ mod 4 \cite[Theorem 7.5]{kirk-klassen},
when the representations $\rho_{\phi_t} = \hol(A_t|_X)$ are
noncentral for all $t$ when restricted to
$T = \partial X$. Recall that $\rho_{\phi_t}|_T$ is central if and only if
$(\alpha_t,\beta_t)\in
(\hz)^2$.

Let us reparametrize $\phi_t$, so that
$\varrho_t:=(\alpha_t,\beta_t)$ lifts to a path $\tilde\varrho_t$ in $\Rt^2$.
The spectral flow on the solid torus in Theorem \ref{splittingformula} has been computed in Theorem
\ref{SFSolidTorus}. Therefore, we focus on the spectral flow on $X$
and the Maslov triple indices in Theorem \ref{splittingformula}.

We apply Proposition \ref{scattidentif} to the
cohomology computations in \cite[Lemma 7.7]{kirk-klassen} to compute the
scattering Lagrangian $\cL_{X,t}$ of $A_t$ on $X$ when $\rho_{\phi_t}|_T$ is
noncentral. Analogous computations yield $\cL_{X,t}$
when $\rho_{\phi_t}|_T$ is central but $\rho_{\phi_t}$ is irreducible. We see that the
scattering Lagrangian on $X$ splits into the $\R i$ and $\C
j$ part $\cL_{X,t}= \hat \cL_{X,t} \oplus \check \cL_{X,t}$, just
like the scattering Lagrangian on $S$. We have
$$
\hat \cL_{X,t} = \text{span} \{i (\det(B+I)\, dm -
c\, dl), i\, dm \wedge dl\}
$$
and
$$
\check \cL_{X,t} =
\begin{cases}
0 & \text{if } (\alpha_t,\beta_t) \notin (\hz)^2\\
e^{i(2 \alpha_t m + 2 \beta_t l)} (\C j\, dm \oplus \C j\, dm \wedge dl)
&
\text{if } (\alpha_t,\beta_t) \in (\hz)^2.
\end{cases}
$$
Thus, by the additivity of the Maslov
triple index, the Maslov triple indices in
Theorem \ref{splittingformula} reduce to $\tau
(K^-_{\alpha_0,\beta_0,i},\check\cL_{S,0},\check\cL_{X,0})$ and $\tau
(K^-_{\alpha_1,\beta_1,i},\check\cL_{S,1},\check\cL_{X,1})$.

\begin{lem} \label{CircleMas}  If $(\alpha_{t^\star},\beta_{t^\star})
  \in
  (\hz)^2$ for some $t^\star \in [0,1]$, then for
  small enough $\epsilon>0$ and varying $t\in [-\epsilon, \epsilon]$ we have $$
\Mas(K^-_{(\alpha_{t^\star},\beta_{t^\star}, \theta e^{t i})},\check \cL_{X,A_{t^\star}}) =
\begin{cases}
2 & \text{if } \theta = \pm 1\\
0 & \text{otherwise}\\
\end{cases}
$$
\end{lem}

\begin{proof} By Theorem \ref{Ktheta}
  we observe
$$\dim(K^-_{(\alpha_{t^\star},\beta_{t^\star},\theta)}\cap \check \cL_{X,A_{t^\star}})=
\begin{cases}
2 & \text{if } \theta = \pm 1,\\
0 & \text{otherwise.}
\end{cases}$$
Consider the case $\theta = 1$. To make notation
  simpler assume $(\alpha_{t^\star},\beta_{t^\star})=(0,0)$. Then  $e^{i(2\alpha_{t^\star} m +
  2\beta_{t^\star} l)} = 1$. Notice that $K^-_{(0,0, 1)}$
  and $\check\cL_{X,A_{t^\star}}$ intersect in $$\text{span}\{ j\, dm
  \wedge dl - k\, dm, k \, dm \wedge dl + j \, dm\}.$$
Consider the constant path  $\tilde L :=
  \text{span} \{ k\, dm, j\, dm\wedge dl \}$ and the path $$L_t = \text{span}
  \{j\, dm\wedge dl-k (\cos t \,  dm + \sin t
  \, dl), j - k (\cos t\, dl - \sin t \, dm)\} \subset
  K^-_{(\alpha,\beta, e^{t i})}$$ of $2$--dimensional Lagrangians in the symplectic subspace
  $$\text{span}\{ j, k\, dm, k\, dl, j\, dm\wedge dl\} \subset
  \cH^{0+1+2}(M;\C j),$$ parametrized by $t\in [-
  \epsilon,\epsilon]$. These intersect at $t=0$ in $\text{span}\{j\, dm
  \wedge dl - k\, dm\}$. We compute
\bea
\lefteqn{\Mas(L_t, \tilde L)}\\
& = &
\Mas(\text{span}\{ j\, dm\wedge dl -(e^{J t}k \, dm), j - (e^{J t} k
\, dl))\},
\tilde L)\\
& = & \Mas(\cL_1 * e^{J t}\text{span}\{ j\, dm \wedge dl - k \, dm, j
- k\, dl\} * \cL_2, \tilde L)\\
& = & \Mas( e^{J t}\text{span}\{ j\, dm \wedge dl - k \, dm, j
- k \,  dl\},\tilde L)\\
&= &1,
\eea
where $*$ denotes composition of paths, and the paths
\bea
\cL_1(t) & := & \text{span}\{ e^{-J t} j \, dm\wedge dl- e^{-J \epsilon}
k\, dm, e^{-J t} j
-e^{-J \epsilon} k \, dl\},\\
\cL_2(t) & := & \cL_1(t-\epsilon)
\eea
are parametrized by
$t\in [0,\epsilon]$. Observe that $\Mas(\cL_i,\tilde L) = 0$ since
$k\, dl \perp \tilde L$.

Now for the orthogonal complement in $\cH^{0+1+2}(M;\C j)$ we have $$\Mas(L_t^\perp,\tilde L^\perp) = \Mas(J L_t, J\tilde L) =
\Mas(L_t,\tilde L) = 1.$$ Thus $$\Mas(K^-_{(0,0, e^{t i})},\check \cL_{X,A_{t^\star}}) = \Mas(L_t \oplus L_t^\perp,\tilde L \oplus
\tilde L^\perp)
= 2.$$
A
similar computation proves the case $\theta = -1$.
\end{proof}

This implies together with Lemma \ref{FlowFlatSolid} and the fact that
$K^-_{\alpha,\beta,-\theta} = K^+_{\alpha,\beta,\theta}$, that the
Maslov triple indices are
either $0$ or $4$ each.

Now we analyze the spectral flow on $X$, when
$\varrho = (\alpha_t,\beta_t)$ passes through the half integer lattice. In the following
proposition we compute the spectral flow for its lift
$\tilde\varrho$ in $\Rt^2$ locally around the half integer lattice.

\begin{prop} Let $[a,b]$ be a maximal interval with $\varrho([a,b])=(\alpha,\beta)\in(\hz)^2$. Then for
  small enough $\epsilon>0$ and varying $t\in [-\epsilon-a,
  b+\epsilon]$ the spectral flow $$\SF(D_{A_t}|_X;\APS^-_{\tilde\varrho_{t}}(\hat
  \cL_{S,t}))$$ is twice the algebraic
  intersection number of $\tilde\varrho|_{[a,b]}$ with
  $\{(\alpha,\beta,\pm 1)\}\subset \Rt^2$.
\end{prop}

\begin{proof} Consider $(\alpha_t,\beta_t)\in (\hz)^2$. The computation
$H^1(X;su(2)_{\rho_{\phi_{t}}})\cong \R^3$ is the same as in the proof of \cite[Lemma 7.7]{kirk-klassen}. For
$(\alpha_t,\beta_t)\notin (\hz)^2$  we have by our computation of
$\cL_{X,t}$, Proposition \ref{scattidentif} and the long exact sequence of of the
pair $(X,T)$ that $\im(H^1(X,T;su(2)_{\rho_{\phi_{t}}}) \to
H^1(X;su(2)_{\rho_{\phi_{t}}})) = 0$. Together with the computations from \cite[Lemma 7.7]{kirk-klassen} we get $$\im(H^1(X,T;su(2)_{\rho_{\phi_{t}}}) \to
H^1(X;su(2)_{\rho_{\phi_{t}}})) =
\begin{cases}
0 & \text{if } \alpha_t \in \hz \text{ or }\beta_t \notin \hz\\
2 & \text{otherwise.}
\end{cases}$$
Then by Proposition \ref{sescor} $D_{A_t}$ has non-resonance level 0 on $X$ when
$\alpha_\phi \in \hz$ or $\beta_\phi \notin \hz$.

Recall
that by $|\tr B| \neq 2$ we have $\det(B+I)\neq 0$ and $c \neq 0$. Thus $\hat\cL_{X,t}$ is
transverse to $\hat \cL_{S,t}$ for all
$t$ and $\theta^\star \neq \pm 1$,
because the $1$--forms of $\hat\cL_{X,t}$ make up the tangent space to the path $i \alpha_t \, dm +
i\beta_t \, dl$ by \cite[Lemma 6.3]{kirk-klassen}.

We have $\tilde\varrho(a)=(\alpha,\beta,-\theta^\star)$ and $\tilde\varrho(b) = (\alpha,\beta,\theta^\star)$. Since
$\theta^\star \neq \pm 1$, we can find an
$\epsilon > 0$ with $\beta_t \notin \hz$ for $t\in [-\epsilon+a,a)
  \cup (b,b+\epsilon]$. Then $D_{A_t}$ has non-resonance level 0, and by Proposition \ref{sescor} we get for $t\in [-\epsilon+a,a)
  \cup (b,b+\epsilon]$
$$
\Lambda_{X,t} \cap (P^-_t \oplus  \cL_{S,t}) \cong \hat\cL_{X,t} \cap \cL_{S,t} = 0.
$$
For $t=a$ and $t=b$ we also have non-resonance level 0. Since
$\theta^\star \neq \pm 1$ we get
$$
\Lambda_{X,t} \cap \APS^-_{\alpha,\beta,\pm\theta^\star} \cong
(\hat\cL_{X,t} \oplus \check\cL_{X,t}) \cap (\hat\cL_{S,t} \oplus K^-_{\alpha,\beta,\pm\theta^\star} ) = 0.
$$
This implies that for $\varrho|_{[-\epsilon+a,a]}$ and
$\varrho|_{[b,b+\epsilon]}$ the spectral flow on $X$ vanishes.

Since we have non-resonance level $0$ on $X$ for all
  $t\in [a,b]$ we have $\Lambda_{X,t}^\infty =
  P^+_{a_t} \oplus \cL_{X,t}$. By Proposition
  $\ref{intersectionstretching}$ and since the Maslov index is
  invariant under a homotopy of paths of Lagragians which preserves the dimension of the intersections at the endpoints,
  we get for varying $t\in [a,b]$
\bea
\lefteqn{\Mas(\APS^-_{\tilde\varrho_t}(\hat
  \cL_{S,t^\star}),\Lambda_{X,t})}\\
 & = & \Mas(\APS^-_{\tilde\varrho_{t}}(\hat
  \cL_S), P^+_{a_t} \oplus \cL_{X,t})\\
& = & \Mas(P^-_{a_t},P^+_{a_t})
 + \Mas(K_{\tilde\varrho_{t}},\check\cL_{X,t}) + \Mas(\hat
  \cL_{S,t}, \hat \cL_{X,t})
\eea
The first and last summand vanish, the second is determined by
  Lemma \ref{CircleMas}. This completes the computation.
\end{proof}

The above proposition together with the computation in the proof of
\cite[Lemma 7.7]{kirk-klassen}, that the spectral flow on $X$ picks up 2
mod 4, whenever $\beta_t \in \hz$ and $\alpha_t \notin \hz$, implies
 that $$\SF(D_{A_t}|_S;\APS^+_{\tilde\varrho_{t}}(J\hat
  \cL_{S,t})) + \SF(D_{A_t}|_X;\APS^-_{\tilde\varrho_{t}}(\hat
  \cL_{S,t}))$$ is a multiple of 4. We summarize.

\begin{thm}\label{SFOnTBundlesOverS1} Let $A_t$ be a path of
  $SU(2)$ connections on a torus bundle over $S^1$, where $A_0$ and
  $A_1$ are flat and irreducible. Then $\SF(D_{A_t}) \equiv 0 \bmod 4$.
\end{thm}

\bibliographystyle{gtart}

\end{document}